\documentclass[11pt,a4paper,reqno]{amsart}
\title{Local Tabularity is Decidable for Bi-Intermediate Logics of Trees and of Co-Trees}

\author{Miguel Martins}
\address{\textbf{Miguel Martins:} Departament de Filosofia\\
Facultat de Filosofia\\
Universitat de Barcelona (UB)\\
Carrer Montalegre, 6, 08001 Barcelona, Spain}
\email{miguelplmartins561@gmail.com}

\author{Tommaso Moraschini}
\address{\textbf{Tommaso Moraschini:} Departament de Filosofia\\
Facultat de Filosofia\\
Universitat de Barcelona (UB)\\
Carrer Montalegre, 6, 08001 Barcelona, Spain}
\email{tommaso.moraschini@ub.edu}

\date{}

\usepackage{mathpazo}  
\usepackage{graphicx}            

\usepackage[utf8]{inputenc}     
\usepackage{amssymb, latexsym, stmaryrd, dsfont, amsmath, amsthm, amsfonts, mathrsfs, amsbsy, mathrsfs, mathtools}            
\usepackage[bottom,symbol]{footmisc}
\usepackage{appendix, verbatim}
\usepackage{etoolbox,graphicx,color}
\usepackage{amscd}
\usepackage{tabularx}
\usepackage{epsfig}
\usepackage{url}
\usepackage{color}
    \usepackage[all, knot]{xy}
        \xyoption{arc}
        \xyoption{web}
\usepackage{tikz}
\usepackage{multicol}
\usepackage[justification=centering]{caption}

\allowdisplaybreaks[1]							

\usepackage[margin=1in]{geometry}

\usepackage{sidecap}                   
\usepackage{bm}                          
\usepackage{enumerate}                   
\usepackage{microtype}  

\makeatletter                            
\def\MT@register@subst@font{\MT@exp@one@n\MT@in@clist\font@name\MT@font@list
 \ifMT@inlist@\else\xdef\MT@font@list{\MT@font@list\font@name,}\fi}
\makeatother
\usepackage[pdftex,bookmarks,bookmarksnumbered,linktocpage,  
         colorlinks,linkcolor=blue,citecolor=blue]{hyperref}
\usepackage{memhfixc}

\usepackage{color}

\definecolor{Salmon}{RGB}{250,128,114}
\definecolor{Crimson}{RGB}{220,20,60}
\definecolor{DarkOrange}{RGB}{255,140,0}
\definecolor{Khaki}{RGB}{240,230,140}
\definecolor{GreenYellow}{RGB}{173,255,47}
\definecolor{MediumSeaGreen}{RGB}{60,179,113}
\definecolor{OliveDrab}{RGB}{107,142,35}
\definecolor{LightSeaGreen}{RGB}{32,178,170}
\definecolor{Aquamarine}{RGB}{127,255,212}
\definecolor{SteelBlue}{RGB}{70,130,180}
\definecolor{Navy}{RGB}{0,0,128}
\definecolor{Purple}{RGB}{128,0,128}
\definecolor{Orchid}{RGB}{218,112,214}
\definecolor{Brown}{RGB}{165,42,42}
\definecolor{Chocolate}{RGB}{210,105,30}
\definecolor{SandyBrown}{RGB}{244,164,96}

\usepackage{tikz}                            
\usetikzlibrary{positioning}
\usetikzlibrary{automata}
 
\usepackage{hyperref}
\hypersetup{colorlinks,citecolor= blue,urlcolor=blue}

\newtheorem{Theorem}{Theorem}[section]
\newtheorem{Lemma}[Theorem]{Lemma}
\newtheorem{Proposition}[Theorem]{Proposition}

\newtheorem{Corollary}[Theorem]{Corollary}

\newtheorem{Claim}{\textbf{Claim}}[Theorem]

\theoremstyle{definition}
\newtheorem{law}[Theorem]{Definition}
\newtheorem{exa}[Theorem]{\textbf{Example}}

\theoremstyle{remark}
\newtheorem{Remark}[Theorem]{Remark}

\newtheorem{exer*}[Theorem]{Exercise*}

\makeatletter
\newcommand\niton{\mathrel{\m@th\mathpalette\canc@l\owns}}
\newcommand\canc@l[2]{{\ooalign{$\hfil#1/\mkern1mu\hfil$\crcr$#1#2$}}}
\newcommand{\bicat}[0]{\operatorname{\mathsf{bi-HA}}}

\newcommand{\J}[0]{\mathcal{J}}

\newcommand{\bipc}[0]{\operatorname{\mathsf{bi-IPC}}}
\newcommand{\ipc}[0]{\mathsf{IPC}}
\newcommand{\bg}[0]{\operatorname{\mathsf{bi-GA}}}

\newcommand{\lc}[0]{\operatorname{\mathsf{bi-GD}}}

\newcommand{\X}[0]{\mathcal{X}}
\newcommand{\Y}[0]{\mathcal{Y}}
\newcommand{\V}[0]{\mathsf{V}}

\newcommand{\bivar}[0]{\operatorname{\mathsf{bi-HA}}}

\newcommand{\p}[0]{ClopUp(\X)}
\newcommand{\A}[0]{\mathbf{A}}
\newcommand{\B}[0]{\mathbf{B}}

\newcommand{\F}[0]{\mathfrak{F}}
\newcommand{\M}[0]{\mathfrak{M}}

\newcommand{\C}[0]{\mathfrak{C}}
\newcommand{\down}[0]{{\downarrow}}
\newcommand{\up}[0]{{\uparrow}}
\newcommand{\SSS}{\mathbb{S}}
\newcommand{\HHH}{\mathbb{H}}
\newcommand{\PPP}{\mathbb{P}}
\newcommand{\ZZ}{\mathbb{Z}}

\newcommand{\si}[0]{{\sim}}

\newcommand{\nesi}[0]{{\neg} \, {\sim}}

\newcommand{\balg}[0]{\text{bi-G\"odel algebra}}
\makeatother

\let\leq=\leqslant
\let\nleq=\nleqslant
\let\geq=\geqslant

\newcommand{\bit}{\begin{itemize}}    
\newcommand{\eit}{\end{itemize}}
\newcommand{\ben}{\begin{enumerate}}
\newcommand{\een}{\end{enumerate}}
\newcommand{\benormal}{\ben[\normalfont 1.]}   
\let\enormal\een
\newcommand{\benroman}{\ben[\normalfont (i)]}  
\let\eroman\een
\newcommand{\benbullet}{\ben[\textbullet]}     
\let\ebullet\een
\newcommand{\vs}{\vspace{.3cm}}

\DeclareRobustCommand{\upcc}{%
  \mathord{\vphantom{\uparrow}\text{%
    \ooalign{%
      $\uparrow$\kern-.515em\raisebox{-1.25ex}{\scalebox{0.88}{$\circ$}}\cr
    }%
  }}%
}
\newcommand{\upc}{\ensuremath \raisebox{1.5pt}{$\upcc$}}

\DeclareRobustCommand{\downcc}{%
  \mathord{\vphantom{\downarrow}\text{%
    \ooalign{%
      $\downarrow$\kern-.515em\raisebox{1.35ex}{\scalebox{0.88}{$\circ$}}\cr
    }%
  }}%
}
\newcommand{\downc}{\ensuremath \raisebox{-1.5pt}{$\downcc$}}

\newcommand{\precc}{{\prec}}

\newcommand{\qq}{^}

\begin{document}

\maketitle

\begin{abstract}
A bi-Heyting algebra validates the G\"odel-Dummett axiom $(p\to q)\vee (q\to p)$ iff the poset of its prime filters is a disjoint union of co-trees (i.e., order duals of trees).\ Bi-Heyting algebras of this kind are called \textit{bi-G\"odel algebras} and form a variety that algebraizes the extension $\lc$ of bi-intuitionistic logic axiomatized by the G\"odel-Dummett axiom.\ In this paper we establish the decidability of the problem of determining if a finitely axiomatizable extension of $\lc$ is locally tabular. 

Notably, if $L$ is an extension of $\lc$, then $L$ is locally tabular iff $L$ is not contained in $Log(FC)$, the logic of a particular family of finite co-trees, called the \textit{finite combs}.
We prove that $Log(FC)$ is finitely axiomatizable. Since this logic also has the finite model property, it is therefore decidable. 
Thus, the above characterization of local tabularity ensures the decidability of the aforementioned problem.

\end{abstract}

\section{Introduction}

The \textit{bi-intuitionistic propositional calculus} $\bipc$ is the conservative extension of the \textit{intuitionistic propositional calculus} $\ipc$ obtained by introducing $\gets$ to the language, a binary connective which behaves dually to $\to$ and is known as \textit{co-implication} (also called \textit{exclusion} or \textit{subtraction}).
To get an intuition for the behavior of the co-implication, we can utilize the Kripke semantics of $\bipc$ \cite{Rauszer6}: if $x$ is a point in a Kripke model $\M$, then for formulas $\phi$ and $\psi$ we have
\[
\M, x \models \phi \gets \psi \text{ iff } \exists y \leq x \; (\M,y \models \phi \text{ and } \M,y\not \models \psi).
\]  

The co-implication $\gets$ gives $\bipc$ significantly greater expressive power than that of $\ipc$.
This is witnessed, for instance, in \cite{Wolter}, where Gödel's embedding of $\ipc$ into the modal logic $\mathsf{S4}$ is extended to an embedding of $\bipc$ into tense-$\mathsf{S4}$, and a version of the Blok-Esakia Theorem \cite{Blok76,Esakia76} is proved, showing that the lattice $\Lambda(\bipc)$ of \textit{bi-intermediate logics} (consistent axiomatic extensions\footnote{From now on we will use \textit{extension} as a synonym of \textit{axiomatic extension}.} of $\bipc$) is isomorphic to that of consistent normal tense logics containing $\mathsf{Grz.t}$ (see also \cite{Cleani21,Stronkowski}).

Furthermore, the addition of this new connective endows $\bipc$ with a symmetry which is notably absent in $\ipc$, since now each connective $\land, \to, \bot$ has its dual $\lor, \gets, \top$, respectively.
This greater symmetry is reflected in the fact that $\bipc$ is algebraized in the sense of \cite{BP89} by the variety $\bicat$ of \textit{bi-Heyting algebras} \cite{Rauszer3} (Heyting algebras whose order duals are also Heyting algebras). Consequently, the lattice $\Lambda(\bipc)$ is dually isomorphic to that of nontrivial varieties of bi-Heyting algebras. The latter, in turn, is amenable to the methods of universal algebra and duality theory since the category of bi-Heyting algebras is dually equivalent to that of \textit{bi-Esakia spaces} \cite{Esakia2} (see also \cite{Bezhan1}).

Motivated by their connection with bi-intuitionistic logic, the theory of bi-Heyting algebras was developed by Rauszer and other in a series of papers (see, e.g., \cite{Beazer,Kohler1,Rauszer2,Rauszer3,Rauszer6,Sanka1}).
However, bi-Heyting algebras also arise naturally in other fields of research.
For example, the lattice of open sets of an Alexandrov space is always a bi-Heyting algebra, and so is the lattice of subgraphs of an arbitrary graph (see, e.g., \cite{Taylor}). 
Similarly, every quantum system can be associated with a complete bi-Heyting algebra \cite{Doring}. 
Many other examples can be found, especially in the field of topos theory \cite{Lawvere1,Lawvere2,Reyes}. 

The thorough investigation of the lattice of \textit{intermediate logics} (consistent extensions of $\ipc$) was a very fruitful topic in nonclassical logic (see, e.g., \cite{Zakha}), but currently $\Lambda(\bipc)$ lacks such an in-depth analysis (for some recent developments in the study of $\bipc$, see, e.g.,  \cite{Badia,BJib,Gore1,Gore2,Shramko}).
In \cite{Paper1}, the authors contributed to the investigation of $\Lambda(\bipc)$ by focusing on the sublattice $\Lambda(\lc)$ of consistent extensions of the \textit{bi-intuitionistic Gödel-Dummett logic}
\[
\lc\coloneqq \bipc + (p\to q) \lor (q \to p).
\]

The formula $(p\to q) \lor (q \to p)$ is called the 
\textit{prelinearity axiom} (or the \textit{Gödel-Dummett axiom}) and over $\ipc$ it axiomatizes the \textit{intuitionistic linear calculus} $\mathsf{LC}$ (or the intuitionistic \textit{Gödel-Dummett logic}).
$\mathsf{LC}$ has been widely studied (see, e.g., \cite{Dummett,Goedel,Horn1,Horn}), and is well-known to be both the intermediate logic of chains (in the sense that it is Kripke complete with respect to the class of \textit{chains}, i.e., linearly ordered Kripke frames) and the intermediate logic of \textit{co-trees} (Kripke frames with a greatest element and whose principal upsets are linearly ordered).
In contrast, while $\lc$ is the bi-intermediate logic of co-trees, it is proved in \cite[Thm. 3.10]{Paper1} that the bi-intermediate logic of chains is its proper extension 
\[
\operatorname{\mathsf{bi-LC}} \coloneqq \bipc + (p \to q) \lor (q \to p) + \neg[ (q\gets p) \land (p \gets q)],
\]
there called the \textit{bi-intuitionistic linear calculus} (see also \cite[Thm. 4.25]{Paper1} for a different axiomatization of $\operatorname{\mathsf{bi-LC}}$).
This suggests that the language of $\bipc$ is more appropriate to study tree-like structures than that of $\ipc$, since $\bipc$ is capable of distinguishing the class of chains from that of co-trees, while $\ipc$ cannot. 
Yet another example of this is that $\ipc$ is well known to be both the intermediate logic of Kripke frames and that of \textit{trees} (order duals of co-trees), and while $\bipc$ is the bi-intermediate logic of Kripke frames \cite{Rauszer6}, it is also shown in \cite[Thm. 3.10]{Paper1} that the bi-intermediate logic of trees is
\[
\operatorname{\mathsf{bi-GD}}\qq \partial \coloneqq \bipc + \neg[ (q\gets p) \land (p \gets q)].
\]

Notably, because of the symmetric nature of bi-intuitionistic logic, all of the results in \cite{Paper1} and in this current paper about extensions of the bi-intermediate logic of co-trees $\lc$ can be rephrased in a straightforward manner as results on extensions of the bi-intermediate logic of trees $\operatorname{\mathsf{bi-GD}}\qq \partial$, by replacing every occurring formula $\varphi$ by its dual $\lnot \varphi^\partial$ ($\varphi^\partial$ is the formula obtained from $\varphi$ by replacing each occurrence of $\land, \lor, \alpha\to \beta, \alpha\gets \beta, \bot, \top$ in $\varphi$ by $\lor, \land, \beta \gets \alpha, \beta \to \alpha, \top, \bot$ respectively) and every algebra or Kripke frame by its order dual. 

There are also other properties of $\lc$ that diverge significantly from those of its intuitionistic fragment $\mathsf{LC}$.
For example, it is known that the lattice $\Lambda(\mathsf{LC})$ of consistent extensions of $\mathsf{LC}$ is a chain of order type $(\omega+1)^\partial$ (see, e.g., \cite{Zakha}), whereas it is shown in \cite[Thm. 4.16]{Paper1} that the lattice $\Lambda(\lc)$ is not a chain and has the cardinality of the continuum. 
It is also well known that $\mathsf{LC}$ is \textit{locally tabular} \cite{Horn} (a logic $L$ is said to be \textit{locally tabular} when $L$ contains, up to logical equivalence, only finitely many formulas in each finite number of variables), but it is an immediate consequence of \cite[Cor. 5.31]{Paper1} that $\lc$ is not.

The main result of \cite{Paper1} was a criterion for local tabularity in extensions of $\lc$.
In the setting of $\Lambda(\lc)$ this notion is intrinsically connected with $FC\coloneqq \{\C_n \colon n \text{ is a positive integer}\}$, the family of the \textit{finite combs} (a particular type of co-trees depicted in Figure \ref{Fig:finite-combs3}),
\begin{figure}[h]
\begin{tikzpicture}
    \tikzstyle{point} = [shape=circle, thick, draw=black, fill=black , scale=0.35]
    \node [label=right:{$x_1'$}] (1') at (1,0) [point] {};
    \node [label=left:{$x_1$}] (1) at (0.5,0.5) [point] {};
    \node [label=right:{$x_2'$}] (2') at (1.5,.5) [point] {};
    \node [label=left:{$x_2$}] (2) at (1,1) [point] {};
    \node [label=above:{$x_n$}] (n) at (1.75,1.75) [point] {};
    \node [label=right:{$x_n'$}] (n') at (2.25,1.25) [point] {};
    
    \draw (1)--(2);
    \draw (1')--(1);
    \draw (2')--(2);
    \draw (n')--(n);
    \draw [dotted] (2)--(n);
\end{tikzpicture}
\caption{The $n$-comb $\C_n$.}
\label{Fig:finite-combs3}
\end{figure}
or, equivalently, with the \textit{logic of the finite combs} $Log(FC)$ (the set of formulas in the language of $\bipc$ which are valid in every finite comb).
This connection is made explicit in \cite[Cor.~5.31]{Paper1}, where it is shown that if $L$ is an extension of $\lc$, then
\begin{equation} \label{crit}
L \text{ is locally tabular iff } L \nsubseteq Log(FC). 
\end{equation}
It is an immediate consequence of this equivalence that $Log(FC)$ is the unique \textit{prelocally tabular} extension of $\lc$, that is, $Log(FC)$ is not locally tabular but all of its proper extensions are so.

The proof of the aforementioned criterion heavily relies on the theories of Jankov and subframe formulas for $\lc$, which were developed in \cite[Sec. 4]{Paper1} (for an overview of these formulas and their use in superintuitionistic and modal logics we refer to \cite{Bezhan2} and \cite{Zakha}, respectively). 
Very briefly, to each finite co-tree $\Y$, we can associate its Jankov (resp. subframe) formula, denoted by $\J(\Y)$ (resp. $\beta(\Y)$). Then, given a \textit{bi-Esakia co-forest} $\X$ (a bi-Esakia space whose underlying poset is a disjoint union of co-trees), the validity of $\J(\Y)$ or $\beta(\Y)$ on $\X$ yields restrictions on the poset structure of $\X$ based on the poset structure of $\Y$ (see Lemmas \ref{jankov lemma} \& \ref{subframe lemma} for more details). 

In this current paper, we prove that the problem of determining if a finitely axiomatizable extension of $\lc$ is locally tabular is decidable (Theorem \ref{main thm}).
We do so by presenting a finite axiomatization of the logic of the finite combs $Log(FC)$.
Since this logic has the finite model property (FMP for short) by definition, it follows that $Log(FC)$ is decidable, whence the equivalence in (\ref{crit}) ensures the decidability of the aforementioned problem.

\begin{figure}[h]
\centering
\begin{tabular}{cccc}
\begin{tikzpicture}
    \tikzstyle{point} = [shape=circle, thick, draw=black, fill=black , scale=0.35]

    \node [label=left:{$d$}] (d) at (-.5,-.25) [point] {};
    \node [label=right:{$e$}] (e) at (.5,-.25) [point] {};
    \node [label=left:{$b$}] (b) at (-.5,.5) [point] {};
    \node [label=right:{$c$}] (c) at (.5,.5) [point] {};
    \node [label=above:{$a$}] (a) at (0,1) [point] {};
    \node [label=below:{\Large{$\F_0$}}] at (0,-.75) [] {};

    \draw (d)--(b)--(a)--(c)--(e);
\end{tikzpicture}
\hspace{1cm}
\begin{tikzpicture}
    \tikzstyle{point} = [shape=circle, thick, draw=black, fill=black , scale=0.35]
    \node [label=left:{$a$}] (a) at (0,1) [point] {};
    \node [label=left:{$b$}] (b) at (0,0) [point] {};
    \node [label=left:{$c$}] (c) at (0,-1) [point] {};
    \node [label=below:{\Large{$\F_1$}}] at (0,-1.5) [] {};

    \draw (c)--(b)--(a);
\end{tikzpicture}
\hspace{1cm}
\begin{tikzpicture}
    \tikzstyle{point} = [shape=circle, thick, draw=black, fill=black , scale=0.35]
    \node [label=above:{$d$}] (d) at (0,0) [point] {};
    \node [label=above:{$c$}] (c) at (0.5,0.5) [point] {};
    \node [label=above:{$b$}] (b) at (01,01) [point] {};
    \node [label=above:{$a$}] (a) at (01.5,01.5) [point] {};
    \node [label=right:{$a'$}] (a') at (2,01) [point] {};
    \node [label=below:{\Large{$\F_2$}}] at (1.2,-.5) [] {};

    \draw (d)--(a)--(a');

\end{tikzpicture}
\hspace{1cm}
\begin{tikzpicture}
    \tikzstyle{point} = [shape=circle, thick, draw=black, fill=black , scale=0.35]

    \node [label=left:{$b$}] (b) at (-.75,0) [point] {};
    \node [label=left:{$c$}] (c) at (0,0) [point] {};
    \node [label=left:{$d$}] (d) at (.75,0) [point] {};
    \node [label=above:{$a$}] (a) at (0,1) [point] {};
    \node [label=below:{\Large{$\F_3$}}] at (0,-.5) [] {};

    \draw (b)--(a)--(c);
    \draw (d)--(a);
\end{tikzpicture}

\end{tabular}
\caption{The co-trees $\F_0$, $\F_1$, $\F_2,\text{ and }\F_3$.}
\label{Fig:the co-trees1}
\end{figure}
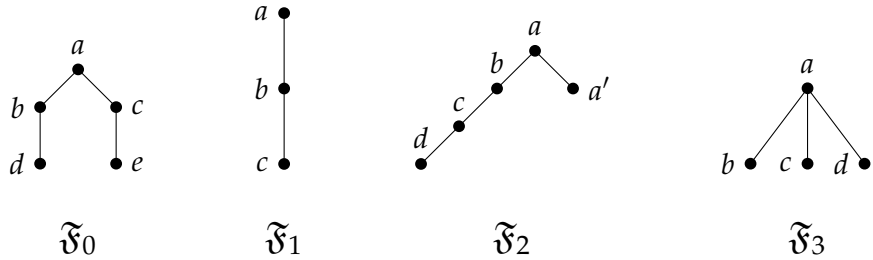

To achieve this finite axiomatization of $Log(FC)$, we take advantage of the rigid and simple structure of the finite combs and of the defining properties of the Jankov and subframe formulas.
We identify four finite co-trees, $\F_0,\F_1,\F_2,\F_3$ (see Figure \ref{Fig:the co-trees1}), such that $\beta(\F_0)$, $\J(\F_1)$, $\J(\F_2)$, and $\J(\F_3)$ are capable of fully capturing the ``comb-like" structure of the bi-Esakia spaces which validate $Log(FC)$, and thus axiomatize this logic.

In more detail, we introduce the logic
\[
    LFC \coloneqq \lc + \beta (\F_0) + \J (\F_1) +\J (\F_2) + \J (\F_3),
    \]
and show that not only does $Log(FC)$ extend $LFC$, by proving that every finite comb validates the axioms of $LFC$ (Proposition \ref{extension of L}), but also that the finite subdirectly irreducible (SI for short) algebras which validate $LFC$ coincide with those which validate $Log(FC)$ (Proposition \ref{prop finite si are equal}).
Section \ref{sec fin gen} is dedicated to the unravelling of the poset structure of dual spaces of finitely generated infinite SI algebras validating $LFC$ (Theorem \ref{Thm poset structure}), and makes crucial use the Coloring Theorem \ref{coloring thm} and of the structure imposed on these spaces by the axioms of $LFC$.
We then use this characterization to prove that $LFC$ enjoys the FMP (Theorem \ref{thm fmp}).
Since so does $Log(FC)$ by definition, and since its finite SI algebraic models coincide with those of $LFC$ by Proposition \ref{prop finite si are equal}, we can conclude that $Log(FC)=LFC$.
It follows that $Log(FC)$ is a finitely axiomatizable logic with the FMP, and is therefore decidable.

\section{Preliminaries}

In this section, we review the basic concepts and results that we will use throughout this paper. For a more in-depth study of $\bipc$ and bi-Heyting algebras, see, e.g., \cite{Paper1,Martins01,Rauszer2,Rauszer3,Rauszer6,Taylor}. 
As a main source for universal algebra we use \cite{Bergman,Sanka2}.

\subsection{Sets and posets}

The class of ordinals will be denoted by $Ord$ and we write $|A|$ for the cardinality of a set $A$.
We will use $\omega$ to denote the first infinite ordinal (i.e., the set of nonnegative integers), while $\mathbb{Z}^+$ denotes the set of positive integers.
The letter $n$ will only be used to represent elements of $\omega$ and the notation $i\leq n$ will always mean either $i\in \{0,\dots ,n\}$ or $i\in \{1,\dots ,n\}$, depending on the context.
We denote the product of a set $A$ with itself (i.e., the set of all pairs of elements of $A$) by $A\qq 2 \coloneqq \{(a,b) \colon a,b \in A\}$.

If $A,B_1,\dots,B_n$ are sets, we sometimes write $A=B_1 \uplus \dots \uplus B_n$ when $A=B_1 \cup \dots \cup B_n$ and the $B_i$'s are pairwise disjoint.
When $R$ is a binary relation on a set $A$ and $\{a,b\} \cup U \subseteq A$, we write $(a,b) \in R$ and $a R b$ interchangeably.
We denote $R[U]\coloneqq \{a\in A \colon \exists u \in U \, (uRa)\}$, and if $U=\{u\}$, we simply write $R(u)\coloneqq \{a \in A \colon u R a\}$. 
One notable binary relation definable on every set $A$ is the \textit{identity relation} $Id_A\coloneqq \{(a,b) \in A\qq 2 \colon a=b \}$.

Let $\X=(X, \leq)$ be a \textit{poset} (i.e., $X$ is a set equipped with a binary relation $\leq$ which is reflexive, transitive, and antisymmetric) and $x,y \in X$.
We write $x<y$ when both $x\leq y$ and $x \neq y$ hold true.
We say that $x$ and $y$ are $\leq$-\textit{comparable} if $x \leq y$ or $y \leq x$.
Otherwise, they are said to be $\leq$-\textit{incomparable}, denoted by $x \perp y$.
When the relation $\leq$ is clear form the context, we shall omit the prefix ``$\leq$", and simply say that $x$ and $y$ are \textit{comparable} (resp. \textit{incomparable}).
If $W,U \subseteq X$ are such that every element of $W$ is incomparable to every element of $U$, we write $W \perp U$, and if $W=\{w\}$, we simply write $w \perp U$.

Let $U \subseteq X$. 
We write $max(U)$ for the set of \textit{maximal} elements of $U$, when viewed as a subposet of $\X$.
If $U$ has a \textit{supremum} in $\X$, we denote it by $sup(U)$.
If, moreover, $sup(U)$ is contained in $U$, we call it the \textit{maximum} (or \textit{greatest element}) of $U$, and write $MAX(U)$ instead.
Similarly, we define the notations $min(U)$, $inf(U)$, and $MIN(U)$, for the set of \textit{minimal} elements of $U$ (when viewed as a subposet of $\X$), for the \textit{infimum} of $U$ in $\X$, and for the \textit{minimum} (or \textit{least element}) of $U$, respectively.

If, when restricted to $U$, the relation $\leq$ is a \textit{linear order} (i.e., any two points of $U$ are comparable), then $U$ is called a $\leq$\textit{-chain}, or simply a \textit{chain} when the relation $\leq$ is clear from the context.
Moreover, if both $MAX(U)$ and $MIN(U)$ exist, then $U$ is a \textit{bounded chain}.
If, for some $n \in \mathbb{Z}^+$, a chain $U$ contains exactly $n$ many elements, we call $U$ an \textit{$n$-chain}.

A map between posets is called an \textit{order embedding} if it is \textit{order invariant} (i.e., it is both order preserving and order reflecting).
It is easy to see that order embeddings are necessarily injective.
An \textit{order isomorphism} is a surjective order embedding.

The poset $\X=(X, \leq)$ is said to be \textit{well-ordered}, or a \textit{well-order}, when $\leq$ is a linear order and every nonempty subset of $X$ has a minimum. 
This is equivalent to $(X,<)$ being order isomorphic to an ordinal equipped with the relation $\in$.

If $x,y\in X$, we say that $x$ is an \textit{immediate predecessor} of $y$, denoted by $x\prec y$, if $x < y$ and no point of $\X$ lies between them (i.e., if $z\in X$ is such that $x \leq z \leq y$, then either $x=z$ or $y=z$). 
If this is the case, we call $y$ an \textit{immediate successor} of $x$. 
If $x$ has a unique immediate successor in $\X$, we denote it by $x^+$. 
We denote the set of points which have immediate successors in $U$ by 
\[
\precc U \coloneqq \{x \in X \colon \exists u \in U \; (x \prec u) \},
\]
and if $U=\{u\}$, we simply write $\precc u \coloneqq \{x \in X \colon x \precc u \}$.

When we write an expression of the form 
\[
(x_1 R_1 x_2 \dots R_{n-1} x_n) \in U,
\]
where $n \in \omega$ and $R_1, \dots R_{n-1}$ are binary relations on $X$, we mean that both 
\[
x_1 R_1 x_2 \dots R_{n-1} x_n \text{ and } x_1, \dots , x_n \in U
\]
hold.
On the other hand, when we write $x_1 R_1 x_2 \dots R_{n-1} x_n \in U$, the only point that is required to be contained in $U$ is $x_n$.

We denote the \textit{upset generated} by $U$ by 
\[
\up U\coloneqq\{x\in X : \exists u\in U \; (u \leq x)\},
\]
and if $U = \up U$, then $U$ is called an \textit{upset}. 
If $U=\{u\}$, we simply write ${\uparrow} u$ and call it a \textit{principal upset}.
The set of upsets of $\X$ will be denoted by $Up(\X)$.
We shall often make use of the following notation
\[
\upc u \coloneqq (\up u) \smallsetminus \{u\} = \{x \in X \colon u < x\}.
\]
The notion of a \textit{downset} and the arrow operators $\down$ and $\downc$ are defined analogously. 
We will always use the convention that the arrow operators and $\precc$ defined above bind stronger than the other set theoretic operations. For example, the expressions $\up U \smallsetminus V$ and $\downc x \cap V$ and $\precc U \smallsetminus V$ are to be read as $(\up U) \smallsetminus V$ and $(\downc x) \cap V$ and $(\precc U) \smallsetminus V$, respectively. 

Finally, when $x,y\in X$, subsets of $X$ of the following forms will be referred to as \textit{intervals}:
\begin{multicols}{2}
    \benbullet
    \item $[x,y] \coloneqq \up x \cap \down y = \{u \in X \colon x \leq u \leq y \}$; \vspace{.2cm}

\item $]x,y[\; \coloneqq \upc x \cap \downc y = \{ u \in X \colon x<u<y \};$

    \item $[x,y[\; \coloneqq \up x \cap \downc y = \{u \in X \colon x \leq u < y \}$; \vspace{.2cm}

    \item $]x,y] \coloneqq \upc x \cap \down y = \{u \in X \colon x < u \leq y \}$.
\ebullet
\end{multicols}

\subsection{Bi-intuitionistic propositional logic} \label{sec bi-ipc}

Let $Fm$ be the set of all formulas in the language $\land, \lor, \to, \gets, \top, \bot$, built up from a denumerable set of variables $Prop$.
Given a formula $\phi$, we write $\lnot \phi$ and $\si \phi$ as a shorthand for $\phi \to \bot$ and $\top\gets\phi$, respectively. 
The \textit{bi-intuitionistic propositional calculus} $\bipc$ is the least subset of $Fm$ that contains the \textit{intuitionistic propositional calculus} $\ipc$ and the eight axioms below, and is moreover closed under modus ponens, uniform substitutions, and the \textit{double negation rule} ``from $\phi$ infer $\nesi \phi$". 
\begin{multicols}{2}
\benormal
    \item $p\to\big(q\vee(p\gets q)\big)$,
    \item $(p\gets q)\to \si (p\to q)$,
    \item $\big((p\gets q )\gets r\big)\to (p\gets q\vee r)$,
    \item $ \neg(p\gets q)\to (p\to q)$,
    \item $\big(p\to (q\gets q)\big)\to \neg p$,
    \item $\neg p\to (p\to (q\gets q)$,
    \item $\big((p\to p)\gets q\big)\to\si q$,
    \item $\si q\to \big((p\to p)\gets q\big)$.
\enormal
\end{multicols}
 \noindent A set of formulas $L$ closed under the three inference rules listed above is called a \textit{super-bi-intuitionistic logic} if it contains $\bipc$. 
Let $L$ be a super-bi-intuitionistic logic and $\phi, \varphi \in Fm$.
We say that $\phi$ is a \textit{theorem} of $L$, denoted by $L\vdash \phi$, if $\phi \in L$. Otherwise, write $L\nvdash \phi$.
If $L \vdash \phi \leftrightarrow \varphi$, then $\phi$ and $\varphi$ are said to be $L$-equivalent formulas.
We call $L$ \textit{consistent} if $L\nvdash \bot$ and \textit{inconsistent} otherwise.

Given another super-bi-intuitionistic logic $L'$, we say that $L'$ is an \textit{extension} of $L$ if $L\subseteq L'$. 
Moreover, $L'$ is proper if $L \neq L'$.
We identify the \textit{classical propositional calculus} \textsf{CPC} with the proper extension of $\bipc$ obtained by adding the \textit{law of excluded middle} $p\lor \neg p$ as an axiom (notably, in \textsf{CPC} the co-implication $\gets$ is term-definable by the other connectives, since $ (p \gets q) \leftrightarrow (p \land \neg q) \in \textsf{CPC}$).  
Consistent extensions of $\bipc$ are called \textit{bi-intermediate logics}, and it can be shown that if $L$ is a super-bi-intuitionistic logic, then $L$ is a bi-intermediate logic iff $ L \subseteq \textsf{CPC}$.
When ordered by inclusion, the set of consistent extensions of a bi-intermediate logic $L$ forms a bounded lattice, which we denote by $\Lambda (L)$.
Given a set of formulas $\Sigma$, we denote by $L+\Sigma$ the least (with respect to inclusion) bi-intuitionistic logic containing $L\cup \Sigma$, and call it the \textit{extension of $L$ axiomatized by $\Sigma$}. If $\Sigma=\{\phi_1, \dots , \phi_n\}$ is finite, we simply write $L+\phi_1+\dots + \phi_n$. 

The \textit{implicative degree} of formulas in the language of $\bipc$ is defined recursively as follows, where $p \in Prop$ and $\phi,\psi \in Fm$:
\benbullet
    \item $ipd(\bot)\coloneqq 0$;
    \item $ipd(p)\coloneqq 0$;
    \item $ipd(\phi \land \psi) \coloneqq MAX\{ipd(\phi),ipd(\psi)\}$;
    \item $ipd(\phi \lor \psi) \coloneqq MAX\{ipd(\phi),ipd(\psi)\}$;
    \item $ipd(\phi \to \psi) \coloneqq MAX\{ipd(\phi),ipd(\psi)\}+1$;
    \item $ipd(\phi \gets \psi) \coloneqq MAX\{ipd(\phi),ipd(\psi)\}+1$.
\ebullet

\subsection{Varieties of algebras}

We denote by $\mathbb{H}, \mathbb{S}$, and $\mathbb{P}$ the class operators of closure under homomorphic images, subalgebras, and (direct) products, respectively. 
A variety $\V$ is a class of (similar) algebras closed under $\mathbb{H}, \mathbb{S}$, and $\mathbb{P}$. By Birkhoff's Theorem, varieties coincide with classes of algebras that can be axiomatized by sets of equations (see, e.g., \cite[Thm.\ II.11.9]{Sanka2}). The smallest variety $\mathbb{V}(\mathsf{K})$ containing a class $\mathsf{K}$ of algebras is called the \textit{variety generated by $\mathsf{K}$} and coincides with $\HHH\SSS\PPP(\mathsf{K})$ (see, e.g., \cite[Thm.\ II.9.5]{Sanka2}).

Given an algebra $\A$, we denote its congruence lattice by $Con(\A)$.\ We say that $\A$ is \textit{subdirectly irreducible}, or SI for short, if $Con(\A)$ has a second least element.

Given a class $\mathsf{K}$ of algebras, we denote by $\mathsf{K}^{<\omega}$, $\mathsf{K}_{SI}$, and $\mathsf{K}_{SI}^{<\omega}$ the classes of finite members of $\mathsf{K}$, SI members of $\mathsf{K}$, and SI members of $\mathsf{K}$ which are finite, respectively. It is a consequence of the Subdirect Decomposition Theorem (see, e.g., \cite[Thm.\ II.8.6]{Sanka2}) that if $\mathsf{K}$ is a variety, then $\mathsf{K} = \mathbb{V}(\mathsf{K}_{SI})$.

A variety $\V$ is said to be \textit{locally finite} if its finitely generated members are finite, and to have the \textit{finite model property} (FMP for short) if it is generated by its finite members.
Using the Subdirect Decomposition Theorem, one can easily show that $\V$ enjoys the FMP iff $\V=\mathbb{V}(\mathsf{\V}_{SI}^{<\omega})$.

\subsection{Bi-Heyting algebras}

A \textit{bi-Heyting algebra} $\A=(A,\land,\lor,\to,\gets,0,1)$ is an algebra with a Heyting algebra reduct $\A\qq -$ whose order dual $(\A\qq-)\qq \partial$ is also a Heyting algebra.
In other words, $\A$ can be viewed as both a Heyting and a co-Heyting algebra.
This is equivalent to $\A$ having a bounded distributive lattice reduct and satisfying the \textit{residuation laws}: for every $a,b\in A$, there are elements $a\to b,$ $a\gets b \in A$ such that
\[
\big( c\leq a\to b \iff a\land c\leq b \big) \text{ and } \big( a\gets b\leq c \iff a\leq b\lor c \big),
\]
for every $c\in A$. 
We will often use the abbreviations $\neg a \coloneqq a \to 0$ and $\si a \coloneqq 1 \gets a$. 

We denote the class of bi-Heyting algebras (which is known to be a variety) by $\bivar$.
Let $\V$ be a variety of bi-Heyting algebras.
When order by inclusion, the set of nontrivial subvarieties of $\V$ forms a bounded lattice, which we denote by $\Lambda (\V)$.

A \textit{bi-Heyting homomorphism} $f \colon \A \to \B$ is a map between bi-Heyting algebras that preserves all of the appropriate algebraic operations.
If $f$ is injective, we call it an \textit{embedding}, and say that $\A$ \textit{embeds} into $\B$.
If $f$ is moreover surjective, we say that $\A$ and $\B$ are \textit{isomorphic}, and denote this by $\A \cong \B$.

Let $\A \in \bivar$ and $a_1, \dots , a_n \in A$. 
We denote the least (with respect to inclusion) bi-Heyting subalgebra of $\A$ containing $\{a_1,\dots,a_n\}$ by $\langle a_1, \dots , a_n \rangle$.
If $\A = \langle a_1, \dots , a_n \rangle$, we say that the elements $a_1, \dots , a_n$ \textit{generate} $\A$, and that $\A$ is \textit{$n$-generated}.
$\A$ is said to be \textit{finitely generated} if it is $m$-generated, for some $m \in \omega$.

A \textit{valuation} on a bi-Heyting algebra $\A$ is a bi-Heyting homomorphism $v\colon \mathbf{Fm} \to \A$, where $\mathbf{Fm}$ denotes the \textit{algebra of formulas} of the language of $\bipc$. 
It is clear that any map $v \colon Prop \to \A$ can be extended uniquely to a valuation on $\A$.
Let $\{\phi\} \cup \Sigma$ be a set of formulas.
We say that $\phi$ is \textit{valid} on $\A$, denoted by $\A \models \phi$, if $v(\phi)=1$ for all valuations $v$ on $\A$. 
On the other hand, if $v(\phi) \neq 1$ for some valuation $v$ on $\A$, we say that $\A$ \textit{refutes} $\phi$ (via $v$), and write $\A \not \models \phi$. 
If every formula in $\Sigma$ is valid on $\A$, we write $\A \models \Sigma$. 
Otherwise, write $\A \not \models \Sigma$.
If $\mathsf{K}$ is a class of bi-Heyting algebras such that $\A \models \Sigma$ for all $\A \in \mathsf{K}$, we write $\mathsf{K} \models \Sigma$. Otherwise, write $\mathsf{K} \not \models \Sigma$.
If $\Sigma=\{\phi\}$, we simply write $K\models \phi$ (resp. $K \not \models \phi$).  

Every bi-intermediate logic $L$ gives rise to a variety $\V_L\coloneqq \{\A \in \bivar \colon \A \models L\}$ of bi-Heyting algebras, called the \textit{variety of $L$}. 
Conversely, with every variety $\V$ of bi-Heyting algebras we can associate its logic $L_\V\coloneqq Log(\V)=\{\phi \in Fm  \colon \V\models \phi \}$.
These correspondences are one-to-one and inverses to each other. 
In fact, when viewed as a deductive system, the logic $\bipc$ is algebraized in the sense of \cite{BP89} by the variety $\bicat$. 
It follows that: for every bi-intermediate logic $L$, the lattice $\Lambda(L)$ of consistent extensions of $L$ is dually isomorphic to the lattice $\Lambda(\V_L)$ of nontrivial subvarieties of $\V_L$; and for every variety of bi-Heyting algebras $\V$, the lattice $\Lambda(\V)$ is dually isomorphic to $\Lambda(L_\V)$.

\subsection{Bi-Esakia spaces} \label{sec bi-esa}

\begin{law} \label{def bi-esa}
    An ordered topological space $\mathcal{X}=(X, \tau, \leq)$ is a \textit{bi-Esakia space} if it is both an Esakia and a co-Esakia space, i.e., when $\X$ is a compact ordered topological space which satisfies the following conditions (where $ClopUp(\mathcal{X})$ denotes the set of its clopen upsets):

\benroman

    \item if $U$ is clopen, then both $\down U$ and $\up U$ are clopen;
    
    \item\label{PSA} \textit{Priestley separation axiom} (PSA for short):  $$\forall x,y\in X\;\big( x \nleq y \implies \exists V\in ClopUp(\X)\;(x\in V \text{ and }  y\notin V)\big).$$ 
\eroman

\noindent When $V$ is a clopen upset satisfying the above display, we say that $V$ \textit{separates} $x$ from $y$.
\end{law}

\begin{law} \label{def bi-p-morphism}
Let $\X=(X, \leq)$ and $\Y=(Y,\leq)$ be posets. A map $f\colon X \to Y$ is called a \textit{bi-p-morphism}, denoted by $f\colon \X \to \Y$, if it satisfies the following conditions:
\benbullet
    \item \underline{\textbf{\textit{Order preserving:}}} $\forall x,z\in X \; \big(x \leq z \implies f(x) \leq f(z)\big)$; \vspace{.1cm}
    
    \item \underline{\textbf{\textit{Up:}}} $\forall x\in X, \forall y\in Y\; \big(f(x) \leq y \implies \exists z\in \up x\; (f(z)=y)\big)$; \vspace{.1cm}
    
    \item \underline{\textbf{\textit{Down:}}} $\forall x\in X, \forall y\in Y\; \big(y \leq f(x) \implies\exists z\in \down x\; (f(z)=y)\big)$.
\ebullet

A continuous bi-p-morphism $f \colon \X \to \Y$ between bi-Esakia spaces is called a \textit{bi-Esakia morphism}.
If there exists a surjective bi-Esakia morphism from $\X$ onto $\Y$, we say that $\Y$ is a \textit{bi-Esakia image} of $\X$, and if there exists a bijective bi-Esakia morphism between $\X$ and $\Y$, then $\X$ and $\Y$ are said to be \textit{isomorphic}, denoted by $\X \cong \Y$.
\end{law}

\begin{exa} \label{finite kripke are bi-esa}
Every finite poset can be viewed as a bi-Esakia space, when equipped with the discrete topology. In fact, since (bi-)Esakia spaces are Hausdorff (see Proposition \ref{bi-esa prop}), this is the only way to view a finite poset as a bi-Esakia space. Furthermore, because maps between sets equipped with the discrete topology are always continuous, we can view every bi-p-morphism between finite posets as a bi-Esakia morphism.
\end{exa}

\begin{Proposition} \label{prop bi-p-morph}
If $f\colon \X \to \Y$ is a bi-p-morphism, then the following conditions hold:
\benroman
    \item $f[\up x]= \up f(x)$ and $f[\down x]=\down f(x)$, for all $x \in X$;
    
    \item $f[max(\X)] \subseteq max(\Y)$ and $f[min(\X)] \subseteq min(\Y)$;
    
    \item if both $MAX(\X)$ and $MAX(\Y)$ exist, then $f\big(MAX(\X)\big)=MAX(\Y)$ and $f$ is necessarily surjective.
\eroman
\end{Proposition}
\begin{proof}
Condition (i) follows immediately from the definition of a bi-p-morphism, while the other two are direct consequences of (i).
\end{proof}

It is known that the celebrated Esakia duality \cite{Esakia2,Esakia3} restricts to a duality between the category of bi-Heyting algebras and bi-Heyting homomorphisms, and that of bi-Esakia spaces and bi-Esakia morphisms. Here, we will only recall the contravariant functors which establish this duality.
Given a bi-Heyting algebra $\A$, we denote its \textit{bi-Esakia dual} by $\A_*\coloneqq (A_*, \tau, \subseteq)$, where $A_*$ is the set of prime filters of $\A$ and $\tau$ is the topology generated by the subbasis
\[
\{\varphi(a)\colon a\in A\}\cup \{A_*\smallsetminus\varphi(a) \colon a\in A\},
\]
 where $\varphi(a)\coloneqq\{F\in A_* \colon a\in F\}$. Notably, it can be shown that $ClopUp(\A_*)=\{\varphi(a)\colon a\in A\}$. Furthermore, if $f \colon \A \to \B$ is a bi-Heyting homomorphism, then its dual is the restricted inverse image map $f_*\coloneqq f^{-1}\colon \B_* \to \A_*$. 
Conversely, if $\X$ is a bi-Esakia space, we denote its \textit{bi-Heyting} (or \textit{algebraic}) \textit{dual} by $\mathcal{X}^*\coloneqq (ClopUp(\mathcal{X}),\cup,\cap,\to,\gets,\emptyset, X),$ where the implications are defined by
\begin{align*}
U\to V&\coloneqq X\smallsetminus \down (U\smallsetminus V)=\big\{x\in X \colon \up x \cap U \subseteq V\big\},\\
U\gets V&\coloneqq \up(U\smallsetminus V)=\big\{x\in X \colon \down x \cap U \nsubseteq V\big\},
\end{align*}
for every $U,V \in \p$.
Moreover, if $f \colon \X \to \Y$ is a bi-Esakia morphism, then its dual is the restricted inverse image map $f^*\coloneqq f^{-1}\colon \Y^* \to \X^*$.

One notable consequence of the bi-Esakia duality that we will often use is: if $\A$ and $\B$ are bi-Heyting algebras, then $\B$ embeds into $\A$ iff $\B_*$ is a bi-Esakia image of $\A_*$ (see, e.g., \cite[Prop.\ 2.15]{Paper1}).
\vs

The following result collects some useful properties of bi-Esakia spaces, most of which will be used in the sequel without further reference.
They are all either well-known results for Esakia spaces (whose proofs can be found in \cite{Esakia3}), or their order dual versions, hence we skip their proof.
Throughout this paper, when we refer to condition (viii) below, we will simply say that ``a bi-Esakia space has enough gaps".

\begin{Proposition} \label{bi-esa prop}
The following conditions hold for every bi-Esakia space $\X$:
\benroman
    \item $\X$ is Hausdorff;

    \item the singletons of $\X$ are closed;

    \item $\X$ is 0-dimensional, that is, $\X$ has a clopen basis;
    
    \item if $Z$ is a closed (resp. open) subset of $\X$, then both $\down Z$ and $\up Z$ are closed (resp. open);

    \item if $x \in X$, then both $\down x $ and $\up x$ are closed;

    \item $min(\X)$ is a closed set;
    
    \item if $Z$ is a closed subset of $\X$ and $x\in Z$, then there are $y \in min(Z)$ and $z \in max(Z)$ satisfying $y \leq x \leq z$;

    \item if $(x <y)\in X$, then there exists a gap between $x$ and $y$, that is, there are $a,b \in X$ such that $x \leq a \prec b \leq y$;

    \item if $C \subseteq X$ is a nonempty chain, then it has both an infimum and a supremum in $\X$, and the closure of $C$ is a bounded chain. 
    In particular, nonempty closed chains in $\X$ are bounded.
    
\eroman 
\end{Proposition}

Given a bi-Esakia space $\X$, a bi-Heyting homomorphism $V\colon \mathbf{Fm} \to \X^*$ is called a \textit{valuation} on $\X$, and the pair $\M\coloneqq (\X,V)$ a \textit{bi-Esakia model} (on $\X$).
If $x \in X$ and $\phi$ is a formula, we say that $x$ satisfies $\phi$ in $\M$ if $x \in V(\phi)$, and denote this by $\M,x \models \phi$.
Otherwise, we say that $x$ \textit{refutes $\phi$ (via $V$)} and write $\M,x \not \models \phi$.
If $V(\phi)=X$ for all valuations $V$ on $\X$, we say that $\phi$ is \textit{valid} on $\X$, or that $\X$ \textit{validates} $\phi$, and write $\X \models \phi$. 
On the other hand, if $V(\phi) \neq X$ for some valuation $V$ on $\X$, we say that $\X$ \textit{refutes $\phi$ (via $V$)}, and write $\M \not \models \phi$ and $\X \not \models \phi$. 
If $\Sigma$ is a set of formulas, then $\X$ \textit{validates} $\Sigma$, written $\X \models \Sigma$, if $\X$ validates every formula in $\Sigma$. 
Otherwise, we say that $\X$ \textit{refutes} $\Sigma$, and write $\X \not \models \Sigma$. 
It is easy to see that $\X \models \phi$ iff $\X^* \models \phi$.\vs

To end this brief discussion on bi-Esakia spaces, we review the Coloring Theorem \cite{EsaGri77}, a result that provides a characterization of the finitely generated bi-Heyting algebras using properties of their bi-Esakia duals. To this end, we need to recall the notions of bi-E-partitions and of colorings on bi-Esakia spaces.

\begin{law} \label{def bi-be}
Let $\X$ be a bi-Esakia space and $E$ an equivalence relation on $X$. We say that $E$ is a \textit{bi-E-partition} of $\X$ if it satisfies the following conditions:
\benbullet
    \item \underline{\textbf{\textit{Up:}}} $\forall x,y,w\in X\, \big(xEy \text{ and } w \in \up y \implies \exists v \in \up x \; (vEw) \big)$; \vspace{.1cm}
    
    \item \underline{\textbf{\textit{Down:}}} $\forall x,y,w\in X\, \big(xEy \text{ and } w \in \down y \implies \exists v \in \down x \; (vEw) \big)$; \vspace{.1cm}
    
    \item \underline{\textbf{\textit{Refined:}}} $\forall x,y \in X \, \big( \neg (xEy) \implies \exists U \in \p \, (E[U]=U \text{ and } |U\cap \{x,y\}|=1) \big).$
\ebullet
A subset $U$ of $X$ satisfying $E[U]= U$ is called \textit{$E$-saturated}. 
Using this terminology, the last condition above can be rephrased as ``any two non-$E$-equivalent elements of $\X$ are separated by an $E$-saturated clopen upset".
We call a bi-E-partition $E$ of $\X$ \textit{trivial} if $E=X^2$, and \textit{proper} otherwise.
\end{law}

Let $\X$ be a bi-Esakia space and $p_1,\dots ,p_n$ a finite number of fixed distinct propositional variables.
Given a map $c\colon \{p_1, \dots ,p_n\} \to ClopUp(\X)$, we associate to each point $x \in X$ the sequence $col(x)\coloneqq (i_1, \dots, i_n)$ defined by
\[
i_k\coloneqq \begin{cases} 1 & \text{if } x\in c(p_k), \\ 0 & \text{if } x\notin c(p_k), \end{cases}
\]
for $k\in \{1,\dots ,n\}$. We call $col(x)$ the \textit{color} of $x$ (relative to $p_1,\dots ,p_n$), the map $c$ an \textit{$n$-coloring} of $\X$, and the pair $(\X,c)$ an \textit{$n$-colored} bi-Esakia space.

Now, if $\A$ is a bi-Heyting algebra endowed with some fixed elements $a_1,\dots ,a_n$, then we can think of this structure as a pair $(\A,v)$, where $v\colon \{p_1,\dots ,p_n\} \to A$ is the map defined by $v(p_i)\coloneqq a_i$, for each $i \leq n$.
Defining $c\colon \{p_1,\dots,p_n\}\to ClopUp(\A_*)$ by $c(p_i)\coloneqq \{x\in A_* \colon a_i \in x\}$, for each $i \leq n$, yields an $n$-coloring of the bi-Esakia dual $\A_*$ of $\A$, and thus an $n$-colored bi-Esakia space $(\A_*,c)$.  \vs 

We are now ready to state the Coloring Theorem \cite{EsaGri77} for bi-Heyting algebras (for a proof, see \cite[Thm.\ 2.19]{Paper1}, noting that it uses a straightforward adaptation of the argument used in \cite[Thm.\ 3.1.5]{Bezhan3} to prove the analogous result for Heyting algebras).

\begin{Theorem}[Coloring Theorem] \label{coloring thm}
 Let $\A$ be a bi-Heyting algebra with elements $a_1,\dots ,a_n$, and let $(\X, c)$ be the corresponding $n$-colored bi-Esakia space. Then $\A= \langle a_1,\dots ,a_n \rangle$ iff every proper bi-E-partition of $\X$ identifies points of different color. 
\end{Theorem}

\subsection{Bi-Gödel Algebras and Co-Trees} \label{sec co-trees}

The \textit{bi-intuitionistic Gödel-Dummett logic} is 
denoted by
\[
\lc \coloneqq \bipc + (p\to q)\lor (q\to p).
\]
This is a bi-intermediate logic and when viewed as a deductive system, it is algebraized in the sense of \cite{BP89} by the variety
\[
\bg \coloneqq \V_{\lc} = \{ \A \in \bivar \colon \A \models (p\to q)\lor (q\to p)\},
\]
whose elements are called \textit{bi-G\"odel algebras} \cite{Paper1, Martins01}. Because of this, there exists a dual isomorphism between the lattice $\Lambda(\lc)$ of consistent extensions of $\lc$ and that of nontrivial subvarieties of $\bg$. \vs

A \textit{co-tree} is a poset $\X$ which has a greatest element (called the \textit{co-root} of $\X$) and whose principal upsets are chains.
A \textit{co-forest} is a (possibly empty) disjoint union of co-trees. 
It follows immediately from the definitions that in a co-tree points have at most one immediate successor, and intervals are always chains.
A \textit{bi-Esakia co-tree} (resp. \textit{bi-Esakia co-forest}) is a bi-Esakia space whose underlying poset is a co-tree (resp. co-forest).

Recall that finite posets can always be viewed as bi-Esakia spaces, when equipped with the discrete topology. 
In particular, so can every finite co-tree.




\begin{Theorem}[{\cite[Thms. 3.1 \& 3.6]{Paper1}}] \label{si balgs}  
Let $\A\in \bivar$. Then $\A$ is a $\balg$ iff $\A_*$ is a bi-Esakia co-forest. Moreover, $\A$ is an SI $\balg$ iff $\A_*$ is a bi-Esakia co-tree.
\end{Theorem}




The theories of Jankov and subframe formulas (for an overview of these formulas and their use in superintuitionistic and modal logics we refer to \cite{Bezhan2} and \cite{Zakha}, respectively) of bi-Gödel algebras were developed in \cite{Paper1, Martins01}. 
With every finite and SI bi-Gödel algebra $\A$, we associate its Jankov formula $\J(\A)$ and its subframe formula $\beta(\A)$. 
If $\Y$ is a finite co-tree, we set $\J(\Y)\coloneqq \J(\Y^*)$ and $\beta(\Y) \coloneqq \beta(\Y^*)$.
The following lemmas illustrate how the refutation or the validity of one of the aforementioned formulas on a bi-Esakia co-forest $\X$ can help us get a better understanding of the poset structure of $\X$.

\begin{Lemma}[Jankov Lemma {\cite[Lem. 4.9]{Paper1}}] \label{jankov lemma}
Let $\X$ be a bi-Esakia co-forest and $\Y$ a finite co-tree.
Then $\X \not \models \J(\Y)$ iff  there exists a surjective bi-Esakia morphism $f \colon \X \twoheadrightarrow \Y$.
\end{Lemma}


\begin{Lemma}[Subframe Lemma {\cite[Lem. 4.24]{Paper1}}] \label{subframe lemma}
Let $\X$ be a bi-Esakia co-forest and $\Y$ a finite co-tree.
Then $\X \not \models \beta(\Y)$ iff $\Y$ order embeds into $\X$.
\end{Lemma}

We end this preliminary part of the paper by presenting two types of equivalence relations on bi-Esakia co-trees that are always proper bi-E-partitions.

\begin{Lemma}[{\cite[Lem. 5.29]{Paper1}}] \label{lemma order-iso}
Let $\X$ be a bi-Esakia co-tree and $w,v\in X$ two distinct points with a common immediate successor. If both $\down w$ and $\down v$ are finite, and there exists an order isomorphism $f\colon \down w \to \down v$, then $E$ is a proper bi-E-partition of $\X$, where
\[
E\coloneqq \big\{(x,y) \in X^2 \colon \big( x\in \down w \text{ and } f(x)=y \big) \text{ or } \big( x\in \down v \text{ and } f(y)=x\big)\big\}\cup Id_X.
\]
\end{Lemma}

\begin{law} \label{def iso chain}
    Let $[y,x]$ be a bounded chain in a bi-Esakia co-tree $\X$.
    We say that $[y,x]$ is an \textit{isolated chain} (in $\X$) when ``no branching occurs in $]y,x]$", i.e., when
    \[
    \down x \smallsetminus [y,x] \subseteq \downc y.
    \]
    In other words, if a point in $[y,x]$ has a predecessor, then said predecessor either lies in the chain $]y,x]$ or lies below $y$.
\end{law}

\begin{Lemma}[{\cite[Lem. 5.28]{Paper1}}] \label{iso chain lem}
    Let $\X$ be a bi-Esakia co-tree and $(y < x) \in X$.
    If $[y,x]$ is an isolated chain, then $E$ is a proper bi-E-partition of $\X$, where
    \[
    E \coloneqq [y,x]^2 \cup Id_X.
    \]
\end{Lemma}

\section{A Finite Axiomatization of the Finite Combs}

\begin{law}
    For each positive integer $n$, we define the $n$-\textit{comb} $\C_n\coloneqq (C_n,\leq)$ as the finite (bi-Esakia) co-tree depicted in Figure \ref{Fig:finite-combs2}.
    We denote the set of the \textit{finite combs} by $FC\coloneqq \{\C_n \colon n \in \mathbb{Z}^+\}$, and its logic by $Log(FC) \coloneqq \{ \phi \in Fm  \colon \forall n \in \mathbb{Z}^+ \; (\C_n \models \phi )\}$.
\end{law}

\begin{figure}[h]
\begin{tikzpicture}
    \tikzstyle{point} = [shape=circle, thick, draw=black, fill=black , scale=0.35]
    \node [label=right:{$x_1'$}] (1') at (1,0) [point] {};
    \node [label=left:{$x_1$}] (1) at (0.5,0.5) [point] {};
    \node [label=right:{$x_2'$}] (2') at (1.5,.5) [point] {};
    \node [label=left:{$x_2$}] (2) at (1,1) [point] {};
    \node [label=above:{$x_n$}] (n) at (1.75,1.75) [point] {};
    \node [label=right:{$x_n'$}] (n') at (2.25,1.25) [point] {};
    
    \draw (1)--(2);
    \draw (1')--(1);
    \draw (2')--(2);
    \draw (n')--(n);
    \draw [dotted] (2)--(n);
\end{tikzpicture}
\caption{The $n$-comb $\C_n$.}
\label{Fig:finite-combs2}
\end{figure}

A bi-intermediate logic $L$ is said to be \textit{locally tabular} when, for every $n \in \omega$, there are only finitely many non-$L$-equivalent formulas in the propositional variables $p_1, \dots , p_n$.
It is easy to see that this is equivalent to the variety $\V_L$ being locally finite.
The interest in the finite combs comes from the following criterion for local tabularity in extensions of $\lc$, as well as the subsequent characterization of $Log(FC)$ as the only \textit{prelocally tabular} extension of $\lc$ (i.e., $Log(FC)$ is not locally tabular, but all of its proper extensions are so).
For a proof of these results, see \cite[Thm.~5.1 \& Cor.~5.31]{Paper1}.

\begin{Theorem}\label{Thm:locally-tabular-main}
An extension $L$ of $\lc$ is locally tabular iff $L \nsubseteq Log(FC)$.
Consequently, $Log(FC)$ is the only prelocally tabular extension of $\lc$.
\end{Theorem}

\begin{figure}[h]
\centering
\begin{tabular}{cccc}
\begin{tikzpicture}
    \tikzstyle{point} = [shape=circle, thick, draw=black, fill=black , scale=0.35]

    \node [label=left:{$d$}] (d) at (-.5,-.25) [point] {};
    \node [label=right:{$e$}] (e) at (.5,-.25) [point] {};
    \node [label=left:{$b$}] (b) at (-.5,.5) [point] {};
    \node [label=right:{$c$}] (c) at (.5,.5) [point] {};
    \node [label=above:{$a$}] (a) at (0,1) [point] {};
    \node [label=below:{\Large{$\F_0$}}] at (0,-.75) [] {};

    \draw (d)--(b)--(a)--(c)--(e);
\end{tikzpicture}
\hspace{1cm}
\begin{tikzpicture}
    \tikzstyle{point} = [shape=circle, thick, draw=black, fill=black , scale=0.35]
    \node [label=left:{$a$}] (a) at (0,1) [point] {};
    \node [label=left:{$b$}] (b) at (0,0) [point] {};
    \node [label=left:{$c$}] (c) at (0,-1) [point] {};
    \node [label=below:{\Large{$\F_1$}}] at (0,-1.5) [] {};

    \draw (c)--(b)--(a);
\end{tikzpicture}
\hspace{1cm}
\begin{tikzpicture}
    \tikzstyle{point} = [shape=circle, thick, draw=black, fill=black , scale=0.35]
    \node [label=above:{$d$}] (d) at (0,0) [point] {};
    \node [label=above:{$c$}] (c) at (0.5,0.5) [point] {};
    \node [label=above:{$b$}] (b) at (01,01) [point] {};
    \node [label=above:{$a$}] (a) at (01.5,01.5) [point] {};
    \node [label=right:{$a'$}] (a') at (2,01) [point] {};
    \node [label=below:{\Large{$\F_2$}}] at (1.2,-.5) [] {};

    \draw (d)--(a)--(a');

\end{tikzpicture}
\hspace{1cm}
\begin{tikzpicture}
    \tikzstyle{point} = [shape=circle, thick, draw=black, fill=black , scale=0.35]

    \node [label=left:{$b$}] (b) at (-.75,0) [point] {};
    \node [label=left:{$c$}] (c) at (0,0) [point] {};
    \node [label=left:{$d$}] (d) at (.75,0) [point] {};
    \node [label=above:{$a$}] (a) at (0,1) [point] {};
    \node [label=below:{\Large{$\F_3$}}] at (0,-.5) [] {};

    \draw (b)--(a)--(c);
    \draw (d)--(a);
\end{tikzpicture}

\end{tabular}
\caption{The co-trees $\F_0$, $\F_1$, $\F_2,\text{ and }\F_3$.}
\label{Fig:the co-trees}
\end{figure}
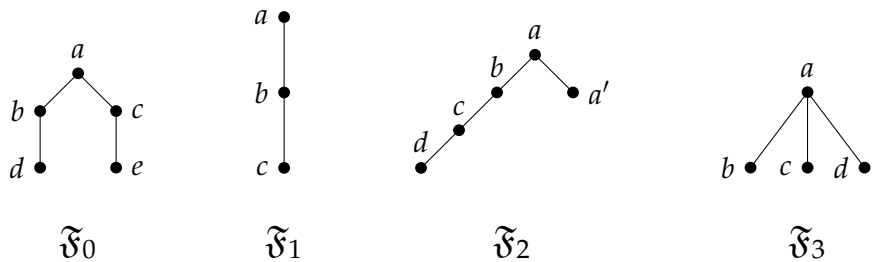

The main goal of this paper is to establish the decidability of the problem of determining when is a
finitely axiomatizable extension of $\lc$ locally tabular.
By the above criterion, this goal can be achieved by showing that $Log(FC)$ is decidable, and since this logic has the FMP by definition, it suffices to prove that it is finitely axiomatizable (because if this is the case, then one can effectively enumerate both the theorems and nontheorems of $Log(FC)$).
Accordingly, this section will be dedicated to proving of the following theorem:

\begin{Theorem} \label{main thm}
    The logic $Log(FC)$ of the finite combs coincides with
    \[
    LFC \coloneqq \lc + \beta (\F_0) + \J (\F_1) +\J (\F_2) + \J (\F_3),
    \]
    where $\F_0, \F_1, \F_2,\F_3$ are the finite co-trees depicted in Figure \ref{Fig:the co-trees}.
    Consequently, $Log(FC)$ is decidable, and so is the problem of determining if a finitely axiomatizable extension of $\lc$ is locally tabular.
\end{Theorem}

Our proof of the above result consists of the following: first we show that not only does $Log(FC)$ extend $LFC$ (Proposition \ref{extension of L}), but also that the finite SI \balg s which validate $LFC$ coincide with those which validate $Log(FC)$ (Proposition \ref{prop finite si are equal}); after unravelling the structure of the dual spaces of the finitely generated SI elements of $\V_{LFC}$ (Theorem \ref{Thm poset structure}), we use this characterization to prove that $LFC$ enjoys the FMP (Theorem \ref{thm fmp}).
Since, by definition, $Log(FC)$ also has the FMP, Proposition \ref{prop finite si are equal} will then entail $Log(FC)=LFC$, as desired.

\subsection{The Validity of the Axioms}

The first step in our proof of Theorem \ref{main thm} is to show $LFC \subseteq Log(FC)$.
We do so by proving that the four axioms of $LFC$ are valid in every finite comb. 
We first recall a helpful lemma.

\begin{Lemma}[{\cite[Lem. 4.15]{Paper1}}] \label{prec}
Let $f\colon \X \to \Y$ be a bi-p-morphism between co-trees. 
Then the co-root of $\X$ must be mapped to the co-root of $\Y$, and if $(x \prec y)\in X$, then either $f(x) = f(y)$ or $f(x) \prec f(y)$.
\end{Lemma} 

\begin{Proposition} \label{extension of L}
    $Log(FC)$ extends the logic $LFC=\lc + \beta (\F_0) + \J (\F_1) +\J (\F_2) + \J (\F_3)$.
\end{Proposition}
\begin{proof}
    By the definitions of these logics, it suffices to show that every finite comb validates the axioms of $LFC$. 
    Let $n$ be a positive integer and consider the $n$-comb $\C_n$.
    It is easy to see that there can be no order embedding from the co-tree $\F_0$ into $\C_n$, hence it follows from the Subframe Lemma \ref{subframe lemma} that $\C_n \models \beta (\F_0)$.

    Next we prove that $\C_n \models \J(\F_1)$. 
    Let us suppose otherwise, so the Jankov Lemma \ref{jankov lemma} forces the existence of a surjective bi-Esakia morphism $f \colon \C_n \twoheadrightarrow \F_1$. 
    By Lemma \ref{prec}, we have $f(x_n')\in \{a,b\}$. 
    But this is already a contradiction: $x_n'$ is a minimal point of $\C_n$, so it must be mapped to a minimal point of $\F_1$ (see Proposition \ref{prop bi-p-morph}.(ii)), but neither $a$ nor $b$ are minimal in $\F_1$.
    Thus, $\C_n$ must validate $\J(\F_1)$.

    To see that $\C_n \models \J(\F_2)$, we again assume otherwise. 
    By the Jankov Lemma \ref{jankov lemma}, there exists a surjective bi-Esakia morphism $f \colon \C_n \twoheadrightarrow \F_2$.
    In particular, this implies $|C_n| \geq 5$, and thus that $n > 2$ by the definition of $\C_n$.
    Using Lemma \ref{prec}, we see that $f(x_n)=a$, hence $f(x_n) \neq b$. 
    As $b$ is not a minimal point of $\F_2$, the aforementioned fact that minimal points must be mapped to minimal points, together with the assumption that $f$ is surjective, implies that the set $\{x_i \in \C_n \colon i < n \text{ and } f(x_i)=b\}$ is nonempty.
    Now, we simply take some $x_i$ belonging to this set, and note that Lemma \ref{prec} forces $f(x_i')\in \{b,c\}$, which is a contradiction: $x_i'$ is minimal, but neither $b$ nor $c$ are so.
    It follows that $\C_n \models \J(\F_2)$.

    Finally, let us suppose that $\C_n \not \models \J(\F_3)$ and arrive at a contradiction.
    By the Jankov Lemma \ref{jankov lemma}, our assumption yields a surjective bi-Esakia morphism $f \colon \C_n \twoheadrightarrow \F_3$. 
    Using Lemma \ref{prec}, we know that $x_n$ must be mapped to $a$.
    It follows that the set $\{x_i \in \C_n \colon f(x_i)=a\}$ is nonempty, and therefore has a least element, by the structure of $\C_n$. 
    Let $x_m$ be this least element.
    Notice that, by combining the aforementioned lemma with Proposition \ref{prop bi-p-morph}.(ii), we can infer that $f(x_m')\in \{b,c,d\}$.
    Moreover, if $m>1$, then the definition of $x_m$ and the structure of $\F_3$ entail that the image of $x_{m-1}$ is a minimal point of $\F_3$.
    As bi-p-morphisms are order preserving, we have $\{f(x_{m-1})\}=f[\down x_{m-1}]$.
    It follows from the structure of $\C_n$, together with $f(x_m)=a$, that
    \[
    |f[\down x_m] \cap \{b,c,d\}|=|f[\{x_m'\} \cup \down x_{m-1}] \cap \{b,c,d\}| \leq 2.
    \]
    Note that $|f[\down x_m] \cap \{b,c,d\}|\leq 2$ also holds when $m=1$, since in this case $x_m'$ is the unique strict predecessor of $x_m$.
    
    But now, for both possibilities on $m$ detailed above, the down condition of bi-p-morphisms (see Definition \ref{def bi-p-morphism}) yields the desired contradiction, since $f(x_m)=a$ implies $\{b,c,d\} \subseteq f[\down x_m]$, and therefore that $|f[\down x_m] \cap \{b,c,d\}| = 3$.
    We conclude $\C_n \models \J(\F_3)$, thus finishing the proof that $LFC \subseteq Log(FC)$.
\end{proof}

\subsection{The Finite SI Members of $\V_{LFC}$}

The next step in our proof of Theorem \ref{main thm} is to show that the finite SI algebras which validate $LFC$ coincide with those which validate $Log(FC)$.

As $LFC \subseteq Log(FC)$ was established above, it is now clear that $\V_{Log(FC)} \subseteq \V_{LFC}$, hence 
\[
(\V_{Log(FC)})_{SI}^{<\omega} \subseteq (\V_{LFC})_{SI}^{<\omega}
\]
follows.
In order to prove the reverse inclusion, we first introduce yet another family of finite bi-Esakia co-trees, and show that in a finite co-tree which validates $\J(\F_3)$, points have at most 2 immediate predecessors.
Recall our notation $\precc x$ for the set of immediate predecessors of a point $x$.

\begin{figure}[h]
\centering
\begin{tabular}{c}
\begin{tikzpicture}
    \tikzstyle{point} = [shape=circle, thick, draw=black, fill=black , scale=0.35]
    \node [label=left:{$y_0$}] (0) at (0,0) [point] {};
    \node [label=right:{$y_1'$}] (1') at (1,0) [point] {};
    \node [label=left:{$y_1$}] (1) at (0.5,0.5) [point] {};
    \node [label=right:{$y_2'$}] (2') at (1.5,.5) [point] {};
    \node [label=left:{$y_2$}] (2) at (1,1) [point] {};
    \node [label=above:{$y_n$}] (n) at (1.75,1.75) [point] {};
    \node [label=right:{$y_n'$}] (n') at (2.25,1.25) [point] {};
    
    \draw (0)--(1)--(2);
    \draw (1')--(1);
    \draw (2')--(2);
    \draw (n')--(n);
    \draw [dotted] (2)--(n);
\end{tikzpicture}
\end{tabular}
\caption{The $n$-hcomb $\C_n'$.}
\label{Fig:finite-hcombs2}
\end{figure}
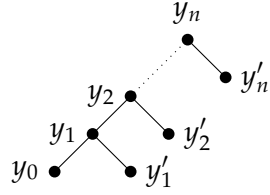

\begin{law} \label{def hcombs}
    For each positive integer $n $, we define the \textit{$n$-comb with handle} (or $n$-\textit{hcomb}, for short) $\C_n'\coloneqq (C_n',\leq)$ as the finite bi-Esakia co-tree depicted in Figure \ref{Fig:finite-hcombs2}.
    Moreover, the $0$-hcomb is defined as the poset with a unique element $y_0$.
    We will refer to these bi-Esakia co-trees as the \textit{finite combs with handle} (or \textit{finite hcombs}, for short).
\end{law}

\begin{Remark} \label{remark combs}
    Notice that, for each $n \in \omega$, there is a surjective bi-p-morphism from the $(n+1)$-comb $\C_{n+1}$ to the $n$-hcomb $\C_n'$, obtained by identifying the points $x_1$ and $x_1'$ of $\C_{n+1}$.
Equivalently, by the bi-Esakia duality, the algebraic dual of $\C_n'$ embeds into the dual of $\C_{n+1}$.
It follows that the algebraic duals of the finite hcombs (which are all SI \balg s, by Theorem \ref{si balgs}) are all contained in $(\V_{Log(FC)})_{SI}^{<\omega}$.
\end{Remark}

\begin{Lemma} \label{lem bounded branching}
    If $\X$ is a finite co-tree which validates $\J(\F_3)$, then $|\precc x| \leq 2$ for all $x \in X$.
\end{Lemma}
\begin{proof}
    We prove the statement by contraposition.
    Suppose that, in a finite co-tree $\X$, there exists a point $x$ which has three distinct immediate predecessors: $z_1, z_2,$ and $z_3$.
    It follows immediately from the assumption that $\X$ is a finite co-tree, together with the definition of immediate predecessors, that 
    \[
    \{\up x, \down z_1, \down z_2, X\smallsetminus (\up x \cup \down z_1 \cup \down z_2)\}
    \]
    is a family of pairwise disjoint nonempty subsets which covers $\X$.
    It is not hard to see that the map $f\colon \X \to \F_3$ defined by $f[\up x ] \coloneqq \{a\}$, $f[\down z_1] \coloneqq \{b\}$, $f[\down z_2] \coloneqq \{c\}$, and $f[X\smallsetminus (\down z_1 \cup \down z_2 \cup \up x)] \coloneqq \{d\}$ is a surjective bi-Esakia morphism.
    By the Jankov Lemma \ref{jankov lemma}, we conclude $\X \not \models \J(\F_3)$.
\end{proof}

Before we proceed, we recall the notion of depth in the setting of finite co-trees.
If $\X$ is a finite co-tree, then we define the \textit{depth} of a point $x \in X$ as $dp(x) \coloneqq | \up x|$, while the \textit{depth} of $\X$ is 
\[
dp(\X) \coloneqq MAX( \{dp(x) \colon x \in X\}) = MAX( \{dp(x) \colon x \in  min(\X) \}).
\]
In other words, the depth of $\X$ is defined as the cardinality of the largest chain contained in $\X$.

\begin{Lemma} \label{lem1}
    If $\A \in (\V_{LFC})_{SI}^{<\omega}$, then $\A_*$ is either isomorphic to a finite comb or to a finite hcomb.
\end{Lemma}
\begin{proof}
    Let $\A \in (\V_{LFC})_{SI}^{<\omega}$.
    We suppose that $\X\coloneqq \A_*$ is not isomorphic to a finite comb, and show that $\X$ must be isomorphic to a finite hcomb.
    If $dp(\X)=1$, then $X$ is a singleton, and is therefore isomorphic to the  0-hcomb $\C_0'$.
    
    Denote the co-root of $\X$ by $y_1$, and suppose that $dp(\X)>1$.
    By our definition of depth, this assumption yields a point $y_2 \in \precc y_1$.
    Observe that if $y_2$ is the sole immediate predecessor of $y_1$, then $y_2$ must be minimal, otherwise we could define a surjective bi-p-morphism from $\X$ onto $\F_1$ (by mapping $y_1 \mapsto a$, $y_2 \mapsto b$, and every other point of $\X$ to $c$). 
    By the Jankov Lemma \ref{jankov lemma}, the existence of such a morphism is equivalent to $\X \not \models \J(\F_1)$, which would contradict the assumption $\A \in \V_{LFC}$, as $\J(\F_1)$ is an axiom of $LFC$.
    But then $\X=(\{y_1, y_2\}, \leq)$ is isomorphic to the $1$-comb, contradicting our other assumption.
    Therefore, $y_1$ must have an immediate predecessor $y_1'$ distinct from $y_2$.
    Since $\J(\F_3)$ is also an axiom of $LFC$, it now follows from Lemma \ref{lem bounded branching} that $\precc y_1 = \{y_1', y_2\}$. 
    If $dp(\X)=2$, then we must have $\precc y_1 = \{y_1', y_2\} \subseteq min(\X)$, and thus that $\X=(\{y_1,y_1',y_2\},\leq)$ is isomorphic to the 1-hcomb $\C_1'$. 
    
    Next we consider the case $dp(\X)>2$, noting that, since $y_1$ is the co-root of $\X$, the definition of depth implies $ \{y_1', y_2\} = \precc y_1 \nsubseteq min(\X)$.
    By symmetry, we assume without loss of generality that $y_2 \notin min(\X)$.

    \begin{Claim}
        For every $y \in X$, if $dp(y)< dp(\X)-1$ then either $y$ is minimal, or $y$ has exactly two immediate predecessors, only one of which is minimal.
    \end{Claim}
    \noindent \textit{Proof of the Claim.} Let $y\in X$ be such that $i\coloneqq dp(y) < dp(\X) -1$.
    We prove the claim by strong induction on $i$. 
    If $i=1$ then $y$ must be the co-root $y_1$, since this is the only point of $\X$ with depth $1$.
    We already established above that $\precc y_1 = \{y_1', y_2\}$ and assumed that $y_2\notin min(\X)$, so it remains to show that $y_1'$ is minimal. 
    But this follows readily from our assumption that $\X$ validates the axiom $\beta(\F_0)$ of $LFC$: if $y_1'$ was not minimal, and since neither is $y_2$, we could define an order embedding from $\F_0$ (see Figure \ref{Fig:the co-trees}) into $\X$, and the Subframe Lemma \ref{subframe lemma} would entail $\X \not \models \beta(\F_0)$. 

    We note that if $dp(\X)=3$, then the base case $i=1$ was the only necessary case to consider.
    Accordingly, we suppose $dp(\X)>3$ and $1 < dp(y)=i < dp(\X)-1$, and assume our induction hypothesis, i.e., that every point in $\X$ whose depth is less than $i$ must be either minimal, or have exactly two immediate predecessors, only one of which is minimal.
    This entails that $\X$ looks like the poset represented in Figure \ref{Fig:fig for claim}, where each $y_{j}'$ is minimal for every $j < i$, but $y_i$ is not.

    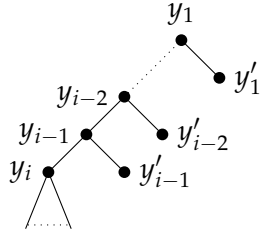
\begin{figure}[h]
\centering
\begin{tabular}{c}
\begin{tikzpicture}
    \tikzstyle{point} = [shape=circle, thick, draw=black, fill=black , scale=0.35]
    \node [label=left:{$y_i$}] (0) at (0,0) [point] {};
    \node [label=right:{$y_{i-1}'$}] (1') at (1,0) [point] {};
    \node [label=left:{$y_{i-1}$}] (1) at (0.5,0.5) [point] {};
    \node [label=right:{$y_{i-2}'$}] (2') at (1.5,.5) [point] {};
    \node [label=left:{$y_{i-2}$}] (2) at (1,1) [point] {};
    \node [label=above:{$y_1$}] (n) at (1.75,1.75) [point] {};
    \node [label=right:{$y_1'$}] (n') at (2.25,1.25) [point] {};
    
    \draw (0)--(1)--(2);
    \draw (1')--(1);
    \draw (2')--(2);
    \draw (n')--(n);
    \draw [dotted] (2)--(n);
    \draw (-.3,-.75)--(0)--(.3,-.75);
    \draw [dotted] (-.27,-.7)--(.27,-.7);
\end{tikzpicture}
\end{tabular}
\caption{The upper part of the poset $\X$.}
\label{Fig:fig for claim}
\end{figure}

As there are exactly two points in $\X$ which have depth $i$, namely, $y_i$ and $y_{i-1}'$, the equality $dp(y)=i$ forces $y$ to be one of the aforementioned points. 
If $y$ is minimal, i.e., if $y=y_{i-1}'$, we are done, so let us suppose that $y$ is not minimal, i.e., that $y= y_i$.
We will now prove that $y_i$ has exactly two immediate predecessors, only one of which is minimal.
Since $y_i$ is not minimal and $\X$ is finite, $y_i$ must have an immediate predecessor, which we will denote by $z$. 
To see that $z$ cannot be the sole element of $\precc y_i$, let us assume otherwise. 
Since $\X$ is a co-tree, $\up z$ is a chain, hence
\[
dp(z)=dp(y_i)+1=i+1<dp(\X).
\]
Using $dp(z) < dp(\X)$ and the known structure of $\X$ (see Figure \ref{Fig:fig for claim}), we can infer from our assumption on $z$ that it cannot be a minimal point.
It is not hard to see that there is a surjective bi-p-morphism $f \colon \X \twoheadrightarrow \F_2$ (see Figure \ref{Fig:the co-trees}), defined by $f[\up y_{i-1}] \coloneqq \{a\}$, $f[min(\X) \smallsetminus \down y_i)] \coloneqq \{a'\}$, $f(y_i) \coloneqq b$, $f(z) \coloneqq c$ and $f[\downc z ] \coloneqq \{d\}$ (recall our notation $\downc z \coloneqq \down z \smallsetminus \{z\}$, and notice that this set in nonempty, as $z$ is not minimal).
Since, by the Jankov Lemma \ref{jankov lemma}, the existence of $f$ is equivalent to $\X$ refuting $\J(\F_2)$, an axiom of $LFC$, we conclude that $y_i$ must have an immediate predecessor $u$ distinct from $z$. 
It now follows from Lemma \ref{lem bounded branching} that $ \precc y_i=\{z, u\}$.

Finally, let us note that $dp(z)=dp(u)=dp(y_i) +1= i+1 <dp(\X)$, so by the known structure of $\X$ (see Figure \ref{Fig:fig for claim}), at least one point in $\precc y_i$ must be nonminimal. 
Denote this point by $y_{i+1}$, and notice that if the other immediate predecessor of $y_i$ distinct from $y_{i+1}$, which we will denote by $y_i'$, was also nonminimal, we could clearly define an order embedding from $\F_0$ (see Figure \ref{Fig:the co-trees}) into $\X$. By the Subframe Lemma \ref{subframe lemma}, this would contradict the assumption $\X \models \beta(\F_0)$.
We conclude that indeed $y_i$ has exactly two immediate predecessors, only one of which is minimal, as desired \qed
\vs

By the previous claim, we now know that $\X$ looks like the poset represented in Figure \ref{Fig:fig for claim2}, where $d\coloneqq dp(\X)$ and $y_{d-1}$ is a nonminimal point of depth $d-1$.

    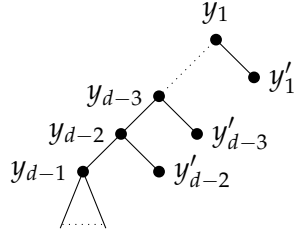
\begin{figure}[h]
\centering
\begin{tabular}{c}
\begin{tikzpicture}
    \tikzstyle{point} = [shape=circle, thick, draw=black, fill=black , scale=0.35]
    \node [label=left:{$y_{d-1}$}] (0) at (0,0) [point] {};
    \node [label=right:{$y_{d-2}'$}] (1') at (1,0) [point] {};
    \node [label=left:{$y_{d-2}$}] (1) at (0.5,0.5) [point] {};
    \node [label=right:{$y_{d-3}'$}] (2') at (1.5,.5) [point] {};
    \node [label=left:{$y_{d-3}$}] (2) at (1,1) [point] {};
    \node [label=above:{$y_1$}] (n) at (1.75,1.75) [point] {};
    \node [label=right:{$y_1'$}] (n') at (2.25,1.25) [point] {};
    
    \draw (0)--(1)--(2);
    \draw (1')--(1);
    \draw (2')--(2);
    \draw (n')--(n);
    \draw [dotted] (2)--(n);
    \draw (-.3,-.75)--(0)--(.3,-.75);
    \draw [dotted] (-.27,-.7)--(.27,-.7);
\end{tikzpicture}
\end{tabular}
\caption{The upper part of the poset $\X$.}
\label{Fig:fig for claim2}
\end{figure}

As $\X$ is a finite co-tree and $y_{d-1}$ is a nonminimal point, it must have immediate predecessors, which will have depth $d=dp(\X)$, and are therefore necessarily minimal. 
Since we assumed that $\X$ is not isomorphic to a finite comb, $y_{d-1}$ must have more than one immediate predecessor (otherwise we would have $\X \cong \C_{d-1}$), so Lemma \ref{lem bounded branching} and the assumption that $\X \models \J(\F_3)$ now imply that $y_{d-1}$ has exactly two such predecessors. 
It is now clear that $\X$ is isomorphic to the finite hcomb $\C_{d-1}'$, as desired.
\end{proof} 

\begin{Proposition} \label{prop finite si are equal}
    The classes $(\V_{Log(FC)})_{SI}^{<\omega}$ and $(\V_{LFC})_{SI}^{<\omega}$ coincide, and their elements are exactly the bi-Heyting algebras whose bi-Esakia duals are either isomorphic to finite combs or to finite hcombs.
\end{Proposition}
\begin{proof}
    As previously discussed, the left to right inclusion is immediate from the fact $LFC \subseteq Log(FC)$ proved in Proposition \ref{extension of L}.
    The reverse inclusion and the second part of the statement follow from the previous lemma, together with Remark \ref{remark combs}.
\end{proof}

With the previous result now proven, we have finished the first part of the proof of Theorem \ref{main thm}.
To conclude that indeed $Log(FC)=LFC$, it now suffices to show that $LFC$ enjoys the FMP, since, by definition, so does $Log(FC) = \{ \phi \in Fm  \colon \forall n \in \mathbb{Z}^+ \; (\C_n \models \phi )\}$.

Our strategy to prove that $LFC$ has the FMP is to characterize the bi-Esakia duals of the finitely generated SI members of $\V_{LFC}$, and subsequently use this characterization to develop a method that, given an arbitrary formula $\varphi$ and an arbitrary infinite finitely generated $\A \in (\V_{LFC})_{SI}$ refuting $\varphi$, extracts a subposet $\Y$ of $\A_*$ that not only refutes $\varphi$, but is crucially isomorphic to a finite comb or to a finite hcomb. 
By the previous proposition, the algebraic dual of $\Y$ must be an element of $(\V_{LFC})_{SI}^{<\omega}$, and we can conclude that $LFC$ enjoys the FMP.

\subsection{The Finitely Generated SI Members of $\V_{LFC}$} \label{sec fin gen}

Throughout this section, we will assume that the reader is familiar with most of the concepts and notations introduced in Section 2, and we will sometimes use them without further reference.
Furthermore, we will work with a fixed but arbitrary infinite finitely generated SI bi-Gödel algebra $\A$ which validates $LFC$.
Let $n \in \mathbb{Z}^+$ be such that $\A$ is $n$-generated, and let $\X \coloneqq \A_*$ be the bi-Esakia dual of $\A$. 
It follows that $\X$ is a bi-Esakia co-tree equipped with an $n$-coloring induced by some distinct propositional variables $p_1, \dots , p_n$.
We denote the co-root of $\X$ by $r$. \vs

Our first goal is to show that $\X$ has a ``comb-like" structure, in the sense that it contains a distinguished minimal point $m$ such that every other point $x$ either lies in $\up m$, or is a minimal point with an immediate successor $x^+$ in $\up m$.
Note that this property of $m$ is shared by the points $x_1'$ in our finite combs (see Figure \ref{Fig:finite-combs2}) and by the points $y_1'$ in our finite hcombs (see Figure \ref{Fig:finite-hcombs2}).
To this end, we need to prove a few auxiliary lemmas.

\begin{Lemma} \label{non-min are comp}
    In $\X$, nonminimal points must be comparable.
\end{Lemma}
\begin{proof}
    Let us suppose, with a view to contradiction, that $x$ and $y$ are nonminimal but incomparable points of $\X$.
    Take some $x' \in \downcc x$ and $y' \in \downcc y$.
    We use the defining properties of co-trees to show $\{x',x\} \perp \{y',y\}$, noting that, by symmetry and by our assumption $x \bot y$, it suffices to prove $x' \perp \{y',y\}$.
    
    Firstly, notice that if $x' \leq y'$ then $x' < y$, so $x,y \in \up x'$ follows.
    As principal upsets in $\X$ are chains, this contradicts our assumption that $x$ and $y$ are incomparable.
    On the other hand, if $y' \leq x'$ then $x', y \in \up y'$, so $x'$ and $y$ must be comparable.
    As proved above, we cannot have $x' \leq y$, and if $y \leq x'$, then $y \leq x$, another contradiction. We conclude that $x' \bot y'$, and note that during our argument, we also showed $x' \bot y$.
    Thus, we have $x' \perp \{y',y\}$, and conclude $\{x',x\} \perp \{y',y\}$.

    Now, since $x \perp y$, it follows that both points are distinct from the co-root $r$ of $\X$.
    Hence, we have $x' < x < r > y > y'$.
    This, together with $\{x',x\} \perp \{y',y\}$, makes clear the existence of an order embedding from $\F_0$ (see Figure \ref{Fig:the co-trees}) into $\X$.
    But this is equivalent to $\X \not \models \beta(\F_0)$ by the Subframe Lemma \ref{subframe lemma}, and thus contradicts our assumption that $\X \models LFC$, since $\beta(\F_0)$ is an axiom of $LFC$.
\end{proof}

\begin{Lemma} \label{lem minimal im pred}
    If $x \in X$, then $x$ can only have at most $2^n$ immediate predecessors which are minimal points of $\X$, i.e.,
    \[
    |\precc x \cap min(\X)| \leq 2^n.
    \]
\end{Lemma}
\begin{proof}
    Let $x \in X$ and note that for distinct $w,v \in \precc x \cap min(\X)$, there is a trivial order isomorphism $f\colon \down w  \to \down v$, since $\down w = \{w\}$ and $ \down v=\{v\}$.
    By Lemma \ref{lemma order-iso}, if $E$ is the smallest equivalence relation on $\X$ that identifies $w$ and $v$, then $E$ is a proper bi-E-partition of $\X$.
    As $\X$ is an $n$-colored bi-Esakia space, it now follows from the Coloring Theorem \ref{coloring thm} that distinct points in $\precc x \cap min(\X)$ must have distinct colors.
    Since there are only $2^n$ possible colors in our fixed coloring of $\X$, the result follows.
\end{proof}

Equipped with the previous lemma, we can now prove that $\X$ has the ``comb-like" structure we described at the beginning of this section.
Recall that if $x$ is a point in a co-tree, then $x$ has at most one immediate successor, and, when it exists, we denote it by $x^+$.

\begin{Proposition} \label{prop diagonal}
    $\X$ has a\textbf{} minimal element $m$ which satisfies $X= \upc m \uplus min(\X)$.
    Moreover, for any other minimal element $y$ distinct from $m$, the point $y^+$ exists and is contained in $\upc m$.
\end{Proposition}
\begin{proof}
    Recall that $\X$ is an infinite co-tree with co-root $r$, so $X \smallsetminus min(\X)$ must be nonempty.
    Moreover, Lemma \ref{non-min are comp} ensures that this set is a chain, thus $u \coloneqq inf(X \smallsetminus min(\X))$ exists by Proposition \ref{bi-esa prop}.(ix).
    If $u\in min(\X)$, we set $m \coloneqq u$, noting that $\upc m = \up m \smallsetminus \{m\}= X \smallsetminus min(\X)$ clearly holds, since any point in $\upc m$ is nonminimal and therefore contained in $X \smallsetminus min(\X)$, while $m=u=inf(X \smallsetminus min(\X)) \in min(\X)$ forces every point in $X \smallsetminus min(\X)$ to lie in $\upc m$.
    On the other hand, if $u \not \in min(\X)$ then $\downc u = \down u \smallsetminus \{u\} \neq \emptyset$ and $u = inf(X \smallsetminus min(\X))=MIN(X \smallsetminus min(\X))$.
    It is easy to see that this equality implies 
    \[
    \emptyset \neq \downc u = \precc u \subseteq min(\X).
    \]
    We take some $m\in \precc u$, noting that since $\X$ is a co-tree, we have
    \[
    \upc m = \up u = \up MIN(X \smallsetminus min(\X))= X \smallsetminus min(\X).
    \]
    In both possibilities for $u$ detailed above, we found an $m \in min(\X)$ satisfying $\upc m = X \smallsetminus min(\X)$, thus proving $X= \upc m \uplus min(\X)$.

    It remains to show the last part of the statement.
    Accordingly, suppose that $y \in min(\X) \smallsetminus \{m\}$ and notice that both $\up m$ and $\up y$ are closed chains which contain $r$, because $\X$ is a bi-Esakia co-tree.
    It follows that $\up m \cap \up y$ is also a nonempty closed chain, and therefore $x \coloneqq  MIN(\up m \cap \up y)$ exists by Proposition \ref{bi-esa prop}.(ix).
    Since $x$ lies above the distinct points $m$ and $y$, we have $x \not \in min(\X)$. 
    This not only implies $y<x$, because $y\in min(\X)$, but also that $x \in X \smallsetminus min(\X)= \upc m$, by above.
    We claim that $x= y\qq +$, i.e., that $x$ is the unique immediate successor of $y$.
    For suppose there is a point $z$ such that $y<z<x$.
    Then also $z \not \in min(\X)$ and $z\in \upc m$.
    This, together with $y<z$, entails $z \in \up m \cap \up y$, and now $z<x$ contradicts our definition of $x = MIN(\up m \cap \up y)$.
    We conclude $x=y\qq +\in \upc m$, as desired.
\end{proof}

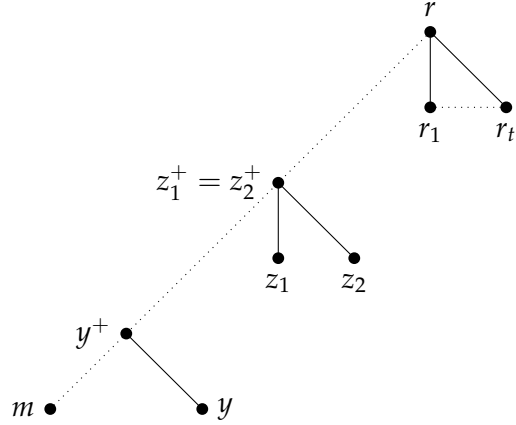
\begin{figure}[h]
\centering
\begin{tabular}{c}
\begin{tikzpicture}
    \tikzstyle{point} = [shape=circle, thick, draw=black, fill=black , scale=0.35]
    \node [label=left:{$m$}] (m) at (0,0) [point] {};
    \node [label=above:{$r$}] (r) at (5,5) [point] {};
    \node [label=below:{$r_1$}] (r1) at (5,4) [point] {};
    \node [label=below:{$r_t$}] (rt) at (6,4) [point] {};
    \node [label=left:{$y^+$}] (y+) at (1,1) [point] {};
    \node [label=right:{$y$}] (y) at (2,0) [point] {};
    \node [label=left:{$z_1^+=z_2^+$}] (z) at (3,3) [point] {};
    \node [label=below:{$z_1$}] (z1) at (3,2) [point] {};
    \node [label=below:{$z_2$}] (z2) at (4,2) [point] {};
    
    \draw[dotted] (m)--(r);
    \draw (y+)--(y);
    \draw (z1)--(z)--(z2);
    \draw (r1)--(r)--(rt);
    \draw[dotted] (r1)--(rt);
\end{tikzpicture}
\end{tabular}
\caption{A rough depiction of the ``comb-like" structure of $\X$.}
\label{Fig:rough depiction}
\end{figure}

\begin{Remark}
    We note that the minimal point $m$ mentioned in the statement of the previous result is not necessarily unique.
    If $u=inf(X \smallsetminus min(\X))\not \in min(\X)$ holds true, then we may have multiple choices for $m \in \precc u$ (although only finitely many by Lemma \ref{lem minimal im pred}). 
    In what follows, we shall work with a fixed $m$, and the chain $\up m$ will sometimes be referred to as ``the diagonal of $\X$".
    Figure \ref{Fig:rough depiction} justifies this terminology.
\end{Remark}

\begin{Corollary} \label{corol diagonal}
    The chain $\upc m$ is an infinite open set.
    Furthermore, for every $x \in \upc m$ there exists a clopen upset $U$ such that $x \in U \subseteq \upc m$ and $U \cap min(\X) = \emptyset$.
\end{Corollary} 
\begin{proof}
As $\X$ is a bi-Esakia space, the set $min(\X)$ is closed by Proposition \ref{bi-esa prop}.(vi).
Because of this, the equality $X= \upc m \uplus min(\X)$ established in Proposition \ref{prop diagonal} entails that $\upc m$ is open.
Moreover, $\upc m$ must be infinite because otherwise, the aforementioned equality and Lemma \ref{lem minimal im pred} would contradict our assumption that $\X$ is infinite.
Since $\X$ is 0-dimensional by Proposition \ref{bi-esa prop}.(iii), it now follows that for every $x \in \upc m$, there exists a clopen $U'$ such that $x \in U' \subseteq \upc m$.
Clearly, $U\coloneqq \up U'$ is a clopen upset satisfying $x \in U \subseteq \upc m$, hence also disjoint from $X \smallsetminus \upc m = min(\X)$, as desired. 
\end{proof}

Our next goal is to prove that the diagonal $\up m$ of $\X$ is well-ordered and to find bounds for the number of immediate predecessors that points in $\up m$ can have. 
Recall that when we write an expression of the form 
\[
(x_1 R_1 x_2 \dots R_{n-1} x_n) \in U,
\]
where $n \in \omega$ and $R_1, \dots R_{n-1}$ are binary relations on $X$, we mean that both $x_1 R_1 x_2 \dots R_{n-1} x_n$ and $x_1, \dots , x_n \in U$ hold.

\begin{Lemma} \label{clop lemma}
    For $x,y,z \in X$, the following conditions hold true:
    \benroman
        \item if $(y \prec x) \in \up m$, then the set $\up x $ is clopen;

        \item if $(x \prec z) \in \upc m$, then the set $\down x$ is clopen;

        \item if $(y \prec x \prec z) \in \up m$, then the set $\{x\}$ is clopen.
    \eroman 
\end{Lemma}
\begin{proof}
    (i) Suppose $(y \prec x) \in \up m$. 
    As $\upc m$ is open by Corollary \ref{corol diagonal} and $\down y$ is closed by Proposition \ref{bi-esa prop}.(v), it follows that $\upc m \smallsetminus \down y = \up x$ is open.
    So the fact that $\up x$ is also closed yields the desired result. \vs

    \noindent (ii) Suppose $(x \prec z) \in \upc m$. 
    As $\upc m$ is open by Corollary \ref{corol diagonal} and $\up z$ is closed by Proposition \ref{bi-esa prop}.(v), we have that $\upc m \smallsetminus \up z = \, ]m,x]$ is an open set.
    It follows from Proposition \ref{bi-esa prop}.(iv) that $\down \, ]m,x]=\down x$ is open, and since we already know that it is closed, we are done.\vs

    \noindent (iii) By the two previous conditions, if $(y \prec x \prec z) \in \up m$ then both $\up x$ and $\down x$ are clopen.
    Therefore, $\up x \cap \down x= \{x\}$ is clopen.
\end{proof}

\begin{law} \label{def up and down limits}
    If $x \in \upc m$ has no immediate predecessors in $\up m$, then we call $x$ an \textit{up limit} (of $\X$).    
If $x \in \up m \smallsetminus \{r\}$ has no immediate successors, then $x$ is called a \textit{down limit} (of $\X$).
\end{law}

\begin{Remark} \label{rem top properties}
    Notice that the above terminology is coherent with certain topological properties of these points.
    For example, if $x$ is an up limit, then it must be a limit point of the set $[y,x]$, for every $y \in [m,x[\,$.
To see this, observe that for such a point $y$, the set $[y,x]=\up y \cap \down x$ is closed and contains $x$, and that $x$ cannot be an isolated point of $[y,x]$, as this would entail that $[y,x[$ is a nonempty closed chain in $\X$.
But in a bi-Esakia space, nonempty closed chains must be bounded (see Proposition \ref{bi-esa prop}), while the definition of up limits ensures that the chain $[y,x[ \,\subseteq \up m$ cannot have a greatest element, since this point would necessarily be an immediate predecessor of $x$ contained in $\up m$.
In a similar way, we can prove that a down limit $x$ must be a limit point of the closed chain $[x,y]$, for every $y \in \, ]x,r]$.
\end{Remark}
 
Our current goal (i.e., showing that the diagonal $\up m$ of $\X$ is well-ordered and finding bounds for the number of immediate predecessors that points in $\up m$ can have) can be achieved by proving the following theorem.
This result, together with Proposition \ref{prop diagonal}, will shed enough light on the structure of the dual spaces of the finitely generated elements of $(\V_{LFC})_{SI}$ for us to prove that $LFC$ has the FMP.
Recall that we are working with a fixed but arbitrary $\A \in (\V_{LFC})_{SI}$, which is both infinite and $n$-generated, and whose dual we denote by $\X$.

\begin{Theorem} \label{thm structure of X}
    The following conditions hold true:
    \benroman
        \item if $x\in \upc m$ is an up limit, then $2 \leq | \precc x | \leq 2^n$;

        \item there are no down limits of $\X$;

        \item $\up m$ is an infinite well-order;

        \item if $x=m^+$, then $| \precc x| \leq 2$;
        
        \item if $x\in \up m \smallsetminus \{m,m^+\}$ and $x$ is not an up limit, then $|\precc x \cap \upc m|=1 = |\precc x \smallsetminus \upc m|$.
    \eroman
\end{Theorem}

We will start by proving condition (i) of the above theorem.
To this end, we will need some technical lemmas, that will also be useful in the sequel.

\begin{Lemma} \label{lem inf chains have teeth}
    If $[y,x] \subseteq \up m$ has more than $2^n$ elements, then there exists $y'\in \, ]y,x]$ such that $\precc y' \smallsetminus \upc m$ is nonempty.
    Consequently, $ (\down x \smallsetminus \down y) \cap min(\X)$ is also nonempty.
\end{Lemma}
\begin{proof}
    Recall that, by Lemma \ref{iso chain lem}, the smallest equivalence relation that identifies the points in an isolated chain (see Definition \ref{def iso chain}) must be a bi-E-partition.
    Since subchains of an isolated chain are clearly isolated chains as well, it follows from the Coloring Theorem \ref{coloring thm} that distinct points of an isolated chain in a colored bi-Esakia space must have distinct colors.
    Consequently, if a bi-Esakia space is $n$-colored, isolated chains have at most $2^n$ elements.

    Now, let us suppose that $[y,x]$ has more than $2^n$ elements and is contained in $\up m$.
    Since $\X$ is $n$-colored, the above discussion entails that $[y,x]$ is not an isolated chain, i.e., that $\down x \smallsetminus [y,x] \nsubseteq \downc y$. 
    Hence, we have $\down x \nsubseteq [y,x] \cup \downc y$, and we can take some $u \in \down x \smallsetminus ([y,x] \cup \downc y)$.
    Notice that $u \notin \up m$, since otherwise $u,y \in \up m$ would entail $u \leq y$ or $y \leq u$, and both conditions clearly contradict our assumption that $u\notin [y,x] \cup \downc y=[y,x] \cup \down y$ (the former would yield $u \in \down y$, while the latter implies $u\in [y,x]$ because we also have $u \leq x$).
    It now follows from Proposition \ref{prop diagonal} that $u$ is necessarily minimal and has a unique immediate successor $u^+$ which is contained in $\upc m$.
    
    Finally, from our assumption $u \in \down x \smallsetminus ([y,x] \cup \downc y)$, we can infer that $u^+ \in \, ]y,x]$.
    To see this, notice that $u < x$ and $u \prec u^+$ entail $u^+ \leq x$, since $\X$ is a co-tree.
    Notice as well that both $y$ and $u^+$ are contained in the chain $\up m$, hence must be comparable points, so $u \nleq y$ and $u \prec u^+$ yield $y < u^+$.
    Since $u\notin \up m$ by above, taking $y'\coloneqq u^+ \in \, ]y,x]$ proves the first part of the statement, because $\precc u^+ \smallsetminus \upc m$ contains $u$.
    Furthermore, Proposition \ref{prop diagonal} ensures $\precc u^+ \smallsetminus \upc m \subseteq min(\X)$, while $y < u^+ \leq x$ and the definition of immediate successors in co-trees imply $\precc u^+ \smallsetminus \upc m \subseteq \down x \smallsetminus \down y$.
    Thus $u \in \precc u^+ \smallsetminus \upc m \subseteq (\down x \smallsetminus \down y) \cap min(\X)$ and we are done.
\end{proof}

Recall that $\X$ is equipped with an $n$-coloring $c\colon \{p_1,\dots,p_n\} \to ClopUp(\X)$ (see Section \ref{sec bi-esa}).
Notice that each element $x \in X$ is contained in the clopen corresponding to $col(x)$, that is, the clopen
\[
C_x \coloneqq \big( \bigcap \{c(p_i) \colon i \leq n \text{ and } x \in c(p_i)\} \big) \smallsetminus \big( \bigcup \{ c(p_i) \colon i \leq n \text{ and } x \notin c(p_i)\} \big).
\]
It is easy to see that not only every point in $C_x$ has the same color, but also that the intersection of this clopen with the diagonal $\up m$ is an interval, as $C_x$ can be written as a finite intersection of upsets and downsets.
It follows that if $(y \leq  z)\in \up m \cap C_x$, then $[y,z] \subseteq C_x$.

Recall as well that every bi-Esakia space has enough gaps (see Proposition \ref{bi-esa prop}).
This fact, and the above discussion on the clopens $C_x$, will be used repeatedly without further reference in what follows.

\begin{Lemma} \label{lem forbidden conditions}
    If $(y \leq x) \in \up m$, then at least one of the following conditions must fail:
    \benroman
        \item $col(y)=col(x)$;

        \item $\precc y \smallsetminus \upc m \neq \emptyset$;

        \item the set $M \coloneqq \big[ (\down x \smallsetminus \down y) \cup (\precc y \smallsetminus \upc m) \big]  \cap min(\X)$ only contains points of the same color and has at least two elements.
    \eroman
\end{Lemma}
\begin{proof}
    Suppose that $(y \leq x) \in \up m$ are such that all of the above conditions hold true.
    We show that this forces
    \[
    E \coloneqq [y,x]^2 \cup M^2 \cup Id_X
    \]
    to be a proper bi-E-partition of $\X$ (see Definition \ref{def bi-be}) that only identifies points of the same color, thus contradicting the Coloring Theorem \ref{coloring thm}.

    We claim that every point in $M$ has an immediate successor contained in $[y,x]$.
    For suppose that $u \in M$.
    If $u=m$, then since $( y \leq x ) \in \up m$ entails $m \notin \down x \smallsetminus \down y$, the definition of $M$ forces $m \in \precc y \smallsetminus \upc m$.
    So $y$ is an immediate successor of $u$ contained in $[y,x]$, satisfying the claim.

    Suppose now that $u \in M$ is distinct from $m$, so in particular we have $u \in min(\X) \smallsetminus \{m\}$.
    It follows from Proposition \ref{prop diagonal} that $u^+$ exists and is contained in $\upc m$.
    By the definition of $M$, either $u \in \precc y \smallsetminus \upc m$ or $u \in \down x \smallsetminus \down y$.
    The former condition yields $u^+=y$ and satisfies our claim, so let us assume the latter, noting that it implies $y < x$, and, in particular, that $x$ is not minimal.
    As $u$ is minimal and lies below $x$, it follows $u <x$.
    By the definition of immediate successors in co-trees, $u^+ \leq x$ follows.
    Furthermore, by the same definition, $u \notin \down y$ and $u < u^+$ entail $u^+ \nleq y$.
    This forces $y < u^+$, since $u^+$ and $y$ must be comparable, as they both lie in the chain $ \up m$. Thus $u^+ \in \,]y,x]$ and we have established our claim.

    Let us now show that $E$ is a proper bi-E-partition which only identifies points of the same color, hence contradicts the Coloring Theorem.
    That this is a proper equivalence relation which only identifies points of the same color is clear from the definition of $E$ and conditions (i) (recall that this condition implies that $[y,x] \subseteq C_x$) and (iii).
    
    To see that $E$ satisfies the up condition, which can be written as
    \[
    \text{ if } u E w \text{ and } v \in \upc u, \text{ then } \exists z \in \up w \text{ s.t. } v E z,
    \]
    let us take $u,w,v$ satisfying the premise of the above display.
    For the sake of nontriviality, we also assume that $(u,v) \notin E$ and $u \neq w$, hence either $u,w\in [y,x]$ or $u,w \in M$, by the definition of $E$.
    If $u$ and $w$ are contained in the chain $[y,x]$, then $(u,v)\notin E$ translates to $v \notin [y,x]$.
    From our assumptions $u\in [y,x]$ and $u<v$, we have $y \leq u < v$, so both $x$ and $v$ are contained in the chain $\up y$.
    Since $\X$ is a co-tree, the points $x$ and $v$ must be comparable.
    But $v \leq x$ would contradict $v \notin [y,x]$, hence $x < v$ must hold true.
    As we also assumed $w \in [y,x]$, we can now take $z \coloneqq v$, since in this case we have $w \leq x < v$ and $vEz$.

    If $u,w \in M$, then by the above claim both $u^+$ and $w^+$ exist and are contained in $[y,x]$, thus $(u^+,w^+)\in E$. 
    By the definition of immediate successors in co-trees, $u <v$ forces $u < u^+ \leq v$.
    So $u^+ \in [y,x]$ now implies that $v$ and $x$ must be comparable.
    If $x \leq v$, then since $w^+\in [y,x]$, we have $w < w^+ \leq v$ and take $z \coloneqq v$.
    If $v \leq x$, then $u^+\in [y,x]$ and $u^+\leq v$ yield $v \in [y,x]$, so $(v,w^+)\in E$ by the definition of $E$, and we take $z \coloneqq w^+$.

    Next we show that $E$ satisfies the down condition, which can be written as
    \[
    \text{ if } u E w \text{ and } v \in \downc u, \text{ then } \exists z \in \down w \text{ s.t. } v E z.
    \]
    Again, without loss of generality, we take distinct $u$ and $w$ such that $u E w$, and $v \in \downc u$ such that $(v,u) \notin E$.
    Moreover, since points in $M$ are minimal by definition, we cannot have $u \in M$ because we assumed that $\downc u \neq \emptyset$.
    Therefore, it suffices to consider the case where $u, w \in [y,x]$.
    As this is a chain and the points $u$ and $w$ were assumed to be distinct, we have either $u < w$ or $w < u$.
    If $u < w$, or, more in general, if $v < w$, we can simply take $z \coloneqq v$ and we are done.
    Suppose now that $w< u$ and $v \nleq w$.
    Notice that these assumptions, together with $(v,u) \notin E$ (which in this case, translates to $v \notin [y,x]$), imply that $v\in (\down u \smallsetminus \down w) \cap min(\X)$, since we have a clear description (see Proposition \ref{prop diagonal} and Figure \ref{Fig:rough depiction}) of $\down u$ as
    \[
    \down u = [w,u] \biguplus \downc w \biguplus \big[ (\down u \smallsetminus \down w) \cap min(\X) \big].
    \]
    As $(w < u) \in [y,x]$, it is clear that $\down u \smallsetminus \down w \subseteq \down x \smallsetminus \down y$, hence 
    \[
    v \in (\down u \smallsetminus \down w) \cap min(\X) \subseteq \big[ (\down x \smallsetminus \down y) \cup (\precc y \smallsetminus \upc m) \big]  \cap min(\X) = M
    \]
    follows.
    Since we assumed that condition (ii) holds, there exists $z \in \precc y \smallsetminus \upc m \subseteq M$.
    It follows from the definition of $E$ that $v E z$, because $v,z\in M$.
    But now, since $w \in [y,x]$, we have $y \leq w$, hence $z < w$ and we are done. 

    It remains to show that $E$ is refined, that is, for
    every $u,v \in \X$,
    \[
    \text{ if } (u,v)\notin E, \text{ then there exists an $E$-saturated clopen upset }  U  \text{ s.t. } |U \cap \{u,v\}|=1.
    \]
    Recall that a clopen $U$ of $\X$ is said to be $E$-saturated when it coincides with 
    \[
    E[U] = \{w\in X \colon \exists z \in U \, (zEw)\}.
    \]
    Accordingly, let $u,v \in X$ be such that $(u,v)\notin E$, and observe that this is equivalent to demanding that 
    \[
    |[y,x] \cap \{u,v\}| \leq 1 \geq |M \cap \{u,v\}| \text{ and } u \neq v.
    \]
    We proceed by cases, recalling that in a bi-Esakia space, there are enough gaps (see Proposition \ref{bi-esa prop}.(viii)) and the upset generated by a clopen is always a clopen (see Definition \ref{def bi-esa}).
    We also recall the definition of the set 
    \[
    M = \big[ (\down x \smallsetminus \down y) \cup (\precc y \smallsetminus \upc m) \big]  \cap min(\X).
    \]
    
    \benbullet
        \item \underline{\textbf{Case:} $u \in M$ and $v=m\notin M$} 

        By the definition of $M$, the assumption $m \notin M$ entails $m \notin \precc y$.
        Since $y \in \upc m$ follows by condition (ii), there must exist $z$ such that $m < z <y$.
        We take a gap $z \leq a \prec b \leq y$, noting that $(a \prec b)\in \upc m$ holds because $m < z$.
        Then Lemma \ref{clop lemma}.(ii) entails that $\down a$ is a clopen set, which clearly contains $m$ and is moreover disjoint from $M$, because $a \in \upc m$ and $a < y$.
        It is easy to see that $U\coloneqq \up \down a$ satisfies the desired conditions, as it is a clopen upset containing $\{v\} \cup [y,x]$ (since $v=m\in \down a $ and $a < y$) and disjoint from $M$ (since so is $\down a$ and $M\subseteq min(\X)$), a set that contains $u$.
        \vspace{.3cm}

        \item \underline{\textbf{Case: }$u \in min(\X)$ and $v \in min(\X) \smallsetminus (M \cup \{m\})$}

        The assumption $v \in min(\X) \smallsetminus (M \cup \{m\})$, together with Proposition \ref{prop diagonal}, implies $v \notin \up m$ and that $v^+$ exists in the chain $\upc m$.
        By the definition of $M$, the assumption also yields both $v \notin \down x \smallsetminus \down y$ and $v \notin \precc y \smallsetminus \upc m$.
        The former condition entails $v^+\notin \, ]y,x]$, while the latter forces $v^+ \neq y$,
        because $v \notin \up m$.
        This establishes $v^+\notin [y,x]$, and since $v^+,y,$ and $x$ are all contained in the chain $\up m$, we have either  $v^+ < y$ or $x < v^+$.
        
        If $v^+ < y$, we take a gap $v^+ \leq a \prec b \leq y$, noting that $(a \prec b)\in \upc m$.
        By Lemma \ref{clop lemma}.(ii), $\down a$ must be a clopen set, which not only contains $v$, but is also clearly disjoint from $M$.
        It follows that $\up \down a$ is a clopen upset containing $\{v\} \cup [y,x]$ (since $v\in \down a$ and $v < v^+ <y$) and disjoint from $ M$ (since so is $\down a$ and $M \subseteq min(\X)$).
        Now, as $u$ and $v$ are distinct minimal points, we have $v \nleq u$, so by the \hyperref[PSA]{PSA} there exists a clopen upset $V$ separating $v$ from $u$.
        Furthermore, since $V$ is an upset containing $v$, the assumption $v \prec v^+ < y$ entails $[y,x] \subseteq V$.
        We take $U \coloneqq \up \down a \cap V$, noting that it clearly satisfies the desired conditions.

        On the other hand, if $x < v^+$ then $v \nleq x$ follows easily from the structure of $\X$, so by two applications of the \hyperref[PSA]{PSA} and taking the intersection of the resulting clopen upsets, there exists a clopen upset $U$ containing $v$ but omitting both $x$ and $u$.
        Because of this, and using the definition of $M$, it is not hard to see that $U \cap ([y,x] \cup M \cup \{u\}) = \emptyset$, hence $U$ is an $E$-saturated clopen upset that separates $v$ from $u$, as desired.
        \vspace{.3cm}

        \item \underline{\textbf{Case:} $u \in \upc m$ and $v \in min(\X)$}

        Simply note that, as both $u$ and $y$ are contained in the chain $\upc m$, we can apply Corollary \ref{corol diagonal} to $u' \coloneqq MIN\{u,y\}$, which yields a clopen upset $U \subseteq \upc m$ that contains $\{u\} \cup [y,x]$ and is disjoint from $min(\X)$, a set that contains $\{v\} \cup M$.     \vspace{.3cm}

        \item \underline{\textbf{Case:} $u,v \in \upc m$} 
        
        Without loss of generality, we assume $u < v$ and proceed by subcases, noting that our assumptions on $u$ and $v$ yield $|[y,x] \cap \{u,v\} | \leq 1$.
        
        If $v<y$, we take a gap $u \leq a \prec b \leq v$.
        By Lemma \ref{clop lemma}.(i), $\up b$ is a clopen upset which moreover contains $\{v\} \cup [y,x]$ and is disjoint from $\{u\} \cup M$.

        If $u < y \leq v$, we take a gap $u \leq a \prec b \leq y$. By Lemma \ref{clop lemma}.(i), $\up b$ is a clopen upset which moreover contains $\{v\} \cup [y,x]$ and is disjoint from $\{u\} \cup M$.

        If $y \leq u \leq x < v$, we take a gap $x \leq a \prec b \leq v$. By Lemma \ref{clop lemma}.(i), $\up b$ is a clopen upset which moreover contains $v$ and is disjoint from $\{u\} \cup [y,x] \cup M$.

        If $x < u < v$, we take a gap $u \leq a \prec b \leq v$. By Lemma \ref{clop lemma}.(i), $\up b$ is a clopen upset which moreover contains $v$ and is disjoint from $\{u\} \cup [y,x] \cup M$.

        In all of the above subcases, $\up b$ satisfies the desired conditions. 
        Furthermore, since all of $u,v,y,x$ are contained in the chain $\up m$, it follows that these subcases cover all possibilities for the points $u$ and $v$.
    \ebullet

    Using the facts $X = \upc m \uplus min(\X)$ and $M \subseteq min(\X)$, it is easy to see that the cases detailed above cover all the possibilities for the points $u$ and $v$, since these facts clearly imply
    \[
    X = \upc m \uplus (M \cup \{m\} ) \uplus [min(\X) \smallsetminus (M \cup \{m\})].
    \]
    We conclude that $E$ is refined, as desired.
\end{proof}

\begin{Lemma} \label{lem up lims have pred}
    If $x \in X$ is an up limit, then it must have immediate predecessors.
\end{Lemma}
\begin{proof}
    Let us suppose that $\precc x = \emptyset$, for some up limit $x \in \upc m$, and arrive at a contradiction. 
    It follows from Proposition \ref{bi-esa prop} that both $\down x$ and $min(\X)$ are closed, hence so is $\down x \cap min(\X)$.
    As $\down x$ is a downset, it is easy to see that $\down x \cap min(\X)=min(\down x)$, so the set $min(\down x)$ is closed.
    Consider an enumeration $ \{x_i\}_{i\in I} = min(\down x)$, noting that clearly $\{x_i\}_{i\in I} \subseteq min(\X)$.
    Now, since we assumed that $x$ is an up limit, it follows that $]m,x[ \, \neq \emptyset$ and that $[y,x]$ is an infinite chain, for every $y \in \, ]m,x[ \,$.
    We take some $y \in \, ]m,x[ $ and using Lemma \ref{lem inf chains have teeth}, we deduce that $min(\down x)$ contains points distinct from $m$, since 
    \[
    \emptyset \neq min(\X) \cap (\down x \smallsetminus \down y) \subseteq min(\down x)
    \]
    and clearly $m \in \down y$.

    For each $i \in I' \coloneqq \{j \in I \colon x_j \neq m\}$, the fact $x_i \in min(\X)\smallsetminus \{m\}$, together with Proposition \ref{prop diagonal}, entails that $x_i^+$ exists and is contained in $\upc m$.
    Moreover, since we assumed $\precc x = \emptyset$, it is clear that $x_i^+ < x$, and so there exists a gap $(x_i^+ \leq a_i \prec b_i < x) \in \upc m$.
    By applying Lemma \ref{clop lemma}.(ii) to $(a_i \prec b_i) \in \upc m$, we see that $\down a_i$ is a clopen set which contains both $x_i$ and $m$.
    It follows that the family $\{\down a_i\}_{i \in I'}$ is an open cover of the closed set $\{x_i\}_{i \in I}= min(\down x)$.
    By compactness, there exists a finite subcover $\{ \down a_{i_1}, \dots , \down a_{i_t}\}$ of $min(\down x)$, but since the $a_{i_1}, \dots , a_{i_t}$ are all contained in the chain $\up m$, there exists some $j \in \{{i_1}, \dots , {i_t}\}$ satisfying
    \[
    min(\down x)=\{x_i\}_{i\in I} \subseteq \bigcup_{l=1}^t \down a_{i_l}=\down a_j.
    \]
    By our definition, $(a_j < x) \in \upc m$, hence the assumption that $x$ is an up limit yields that $[a_j,x]$ is an infinite chain contained in $\up m$.
    By Lemma \ref{lem inf chains have teeth}, $min(\X) \cap (\down x \smallsetminus \down a_j)$ must be nonempty.
    But this set is clearly contained in $min(\down x)$ and disjoint from $\down a_j$, contradicting the inclusion $min(\down x) \subseteq \down a_j$ in the above display.
\end{proof}

\noindent\textit{Proof of Theorem} \ref{thm structure of X}.(i). We will now prove that if $x$ is an up limit, then $  2 \leq | \precc x | \leq 2^n.$
First of all, by the definition of up limits, $x$ has no immediate predecessor in the diagonal $\up m$, i.e., $\precc x \subseteq X \smallsetminus \up m$.
Using the equality $X = \upc m \uplus min(\X)$ established in Proposition \ref{prop diagonal}, it follows that points in $\precc x$ are necessarily minimal, so Lemma \ref{lem minimal im pred} entails that there are at most $2^n$ such points.
This, together with Lemma \ref{lem up lims have pred}, yields $  1 \leq | \precc x | \leq 2^n$, so it remains to show that $x$ cannot have a sole immediate predecessor.

Accordingly, let us suppose that $\precc x= \{x_1\}$, for some $x_1$, and arrive at a contradiction.
Recall that $x_1$ is contained in the clopen $C_{x_1}$ corresponding to $col(x_1)$ and consider the set $min(\down x) \smallsetminus C_{x_1}$.
    It follows from Proposition \ref{bi-esa prop} that both $\down x$ and $min(\X)$ are closed, hence so is $\down x \cap min(\X)$.
    As $\down x$ is a downset, it is easy to see that $\down x \cap min(\X)=min(\down x)$, so the set $min(\down x)$ is closed.
    Hence $min(\down x) \smallsetminus C_{x_1}$ is also closed.
    We claim that $min(\down x) \smallsetminus C_{x_1}$ is contained in $\down a$, for some $a \in [m,x[\,$.
    For suppose that $m$ is the only point in $min(\down x) \smallsetminus C_{x_1}$.
    Then clearly $a \coloneqq m$ would satisfy our claim.
    On the other hand, if $min(\down x) \smallsetminus C_{x_1}$ contains points distinct from $m$, then we can repeat the same argument used in the proof of Lemma \ref{lem up lims have pred}, using gaps and compactness, to find such a point $a \in [m,x[$ whose downset contains $min(\down x) \smallsetminus C_{x_1}$.

    Now, as $x$ is an up limit, we know by Remark \ref{rem top properties} that it must be a limit point of the closed chain $[a,x]$.
    Because of this, the intersection of the clopen $C_x$ with the set $[a,x[$ must be nonempty, since $C_x$ contains $x$.
    Let $z$ be a point in this intersection, and notice that, by the definition of up limits, $[z,x]$ is an infinite chain, that moreover is contained in the clopen $C_x$.
    It now follows from Lemma \ref{lem inf chains have teeth} that we can find a point $y \in \, ]z,x[$ such that $ \precc y \smallsetminus \upc m \neq \emptyset$.
    
    But now, we have $(y < x)\in \up m$, and that $x$ and $y$ satisfy the following conditions: (i) $col(y) = col(x)$ (since $y\in [z,x] \subseteq C_x$); (ii) $ \precc y \smallsetminus \upc m \neq \emptyset$ (by above); and (iii) the set $M \coloneqq \big[ (\down x \smallsetminus \down y) \cup (\precc y \smallsetminus \upc m) \big] \cap min(\X)$ only contains points of the same color (since we established above that $min(\down x) \smallsetminus C_{x_1} \subseteq \down a$, so $a < y$ implies $M \subseteq C_{x_1}$) and has at least two elements ($M$ contains both the nonempty set $\precc y \smallsetminus \upc m$ and $\{x_1\}$, and clearly $y <x$ entails $x_1 \notin \precc y \smallsetminus \upc m$).
    This contradicts Lemma \ref{lem forbidden conditions}, hence we conclude that $\precc x$ cannot be a singleton, as desired. \qed 
    
    \vs

Next we will prove condition (ii) of Theorem \ref{thm structure of X}, that is, we will show that down limits (see Definition \ref{def up and down limits}) do not exist.
We first need to prove a couple of auxiliary results.

\begin{Lemma} \label{lem compact}
    Let $x \in \up m \smallsetminus \{r\}$ be a down limit, $z \in \upc x$, and 
    \[
    M_x \coloneqq \begin{cases}
    \precc x \smallsetminus \upc m & \text{ if } x \neq m, \\
    \{x\} & \text{ if } x = m.
\end{cases}
    \]
    If $C$ is a clopen containing $M_x$, then there exists a point $z' \in \, ]x,z[$ such that $min(\down z' \smallsetminus \down x) \cup M_x \subseteq C$.
    
\end{Lemma}
\begin{proof}
    Let $x$ be a down limit, $z \in \upc x$, and $C$ a clopen set containing $M_x$.
    As $x$ has no immediate successors, $]x,z]$ has no least element, and the equality $min(\down z \smallsetminus \down x) = (\down z \smallsetminus \down x) \cap min(\X)$ follows because $\X$ is a co-tree.
    In particular, $W \coloneqq min(\down z \smallsetminus \down x ) \smallsetminus C \subseteq min(\X)$.

    We claim that $W$ is a closed set.
    Observe that, since $X=\upc m \uplus min(\X)$ by Proposition \ref{prop diagonal}, we have 
    \[
    X \smallsetminus W =  \down x \,  \cup \, ]x,z] \cup (\down r \smallsetminus \down z) \cup C= \upc m \cup min(\down x)\cup (\down r \smallsetminus \down z) \cup C
    \]
    (the above equalities are easily verified by picturing $\X$ and the occurring sets).
    To establish our claim, we show that every point contained in the right most union of the above display has an open neighbourhood disjoint from $W$.
    
    For suppose $u \in X\smallsetminus W$. 
    The cases $u \in \upc m$ and $u \in C$ are clear, since the former set is open by Corollary \ref{corol diagonal}, while the latter is clopen by assumption, and both sets are clearly disjoint from $W$ in view of the above display.
    
    The case $u \in \down r \smallsetminus \down z$ is also easy to see, since this implies $u \notin \down z$, so the \hyperref[PSA]{PSA} yields a clopen upset $U$ separating $u$ from $z$.
    This clopen upset must be clearly disjoint from $\down z$, and in particular, from $W \subseteq \down z$.

    It remains to consider the case $u \in min(\down x)$.
    Notice that if $x=m$, then we have $min(\down x)=\{m\}=M_x$ by the definition of $M_x$, so our assumption on $C$ yields $u \in M_x \subseteq C$, and we already dealt with this possibility above.
    Accordingly, we can assume, without loss of generality, that $x \neq m$, hence $m < x$ holds, and that $u \in min(\down x) \smallsetminus C$.
    In this case, $M_x=\precc x \smallsetminus \upc m \subseteq C$ holds true, so our above assumption that $u \notin C$ yields $u \notin \precc x \smallsetminus \upc m$.
    As we clearly have $u \in min(\down x) \subseteq min(\X)$, from the equality $X = \upc m \uplus min(\X)$ (see Proposition \ref{prop diagonal}) we can now infer $u \not \prec x$.
    Hence, there exists $v$ satisfying $u < v < x$.
    Now, as $v$ is not minimal, $X = \upc m \uplus min(\X)$ implies $ (v < x)  \in \upc m$.
    This yields a gap $(v \leq a \prec b \leq x) \in \upc m$, and it follows from Lemma \ref{clop lemma}.(ii) that $\down a$ is a clopen neighbourhood of $u$ (recall $u < v \leq a$), that moreover is disjoint from $W = min(\down z \smallsetminus \down x) \smallsetminus C$, because $a < x$.
    This finishes our proof of the claim that $W$ is a closed set.

    Consider an enumeration $\{y_i\}_{i \in I}=W$.
    For each $i \in I$, we have $y_i  \nleq x$ by the definition of $W$, so the \hyperref[PSA]{PSA} yields a clopen upset $U_i$ that separates $y_i$ from $x$.
    Again by the definition of $W$, the clopen $C$ omits $y_i$, hence $V_i \coloneqq \up (U_i \smallsetminus C)$ is also a clopen upset that contains $y_i$ but not $x$.
    It follows that $\{V_i\}_{i\in I}$ is a clopen cover of $W=\{y_i\}_{i\in I}$.
    As $W$ is closed by the previous claim, there must exist a finite subcover $\{V_{i_1}, \dots , V_{i_t}\}$ of $W$ by compactness.
    Since, for each $i_j$, the set $V_{i_j}$ is a clopen upset that contains $y_{i_j}$, it follows from 
    \[
    y_{i_j} \in W = min(\down z \smallsetminus \down x) \smallsetminus C \subseteq \down z
    \]
    that $ z \in V_{i_j}$.
    Thus, $ z$ is contained in the clopen $\bigcup_{j=1}^t V_{i_j}$.
    Since we also have $z \in \up m$ (recall that $z \in \upc x \subseteq \up m$), we infer that $\up m \cap \bigcup_{j=1}^t V_{i_j} $ is a nonempty closed chain in a bi-Esakia space, and is therefore bounded by Proposition \ref{bi-esa prop}.
    Let $z_0 \coloneqq MIN(\up m \cap \bigcup_{j=1}^t V_{i_j})$ and notice that this point lies below $z$, by the comment above.
    Furthermore, $z_0$ must be comparable to $x$, since both these points lie in the chain $\up m$. 
    But $z_0 \leq x$ cannot happen, since $\bigcup_{j=1}^t V_{i_j}$ is an union of upsets that omit $x$.
    Hence, we must have $x < z_0 \leq z$.

    Now, as $x$ is assumed to be a down limit, there exists $z'$ such that $ x < z' < z_0$.
    By the definition of $z_0$, the previous inequality yields that $\bigcup_{j=1}^t V_{i_j}$ is an upset that omits $z'$, and is thus disjoint from $\down z'$.
    Since $W \subseteq \bigcup_{j=1}^t V_{i_j}$ by compactness, we now have $\down z' \cap W = \emptyset$.
    Moreover, $z' \in \, ]x,z[$ clearly implies that $min(\down z' \smallsetminus \down x) \subseteq min(\down z \smallsetminus \down x)$.
    Therefore, the fact that $\down z'$ and $W = min(\down z \smallsetminus \down x) \smallsetminus C$ are disjoint now proves $min(\down z' \smallsetminus \down x) \subseteq C$, and the point $z'$ satisfies the conditions of the statement.
\end{proof}

\begin{Proposition} \label{prop bound for down limits}
    If $x$ is a down limit, then $2 \leq | \precc x \smallsetminus \upc m | \leq 2^n.$
\end{Proposition}
\begin{proof}
Let $x$ be a down limit and recall the equality $X = \upc m \uplus min(\X)$ proved in Proposition \ref{prop diagonal}. 
It follows that any point in $\precc x \smallsetminus \upc m$ is necessarily minimal, so Lemma \ref{lem minimal im pred} entails that this set has at most $2^n$ points.

    Now, let us suppose, with a view to contradiction, that $x$ is such that $|\precc x \smallsetminus \upc m| \leq 1$.
    It follows from Remark \ref{rem top properties} that $x$ is a limit point of the closed chain $\up x$.
    Since $C_x$, the clopen corresponding to $col(x)$, is a clopen neighbourhood of $x$, the previous comment yields a point $z \in \upc x \cap C_x$. 
    Again using our assumption that $x$ is a down limit, we see that $[x,z]$ is an infinite chain, that moreover is contained in $C_x$ (recall that $C_x$ can be written as a finite intersection of upsets and downsets, hence its intersection with $\up m$ is always an interval).
    Set
    \[
    M_x \coloneqq \begin{cases}
    \precc x \smallsetminus \upc m  & \text{ if } x \neq m, \\
    \{x\} & \text{ if } x = m,
\end{cases}
\hspace{1cm} \text{ and } \hspace{1cm}
C\coloneqq \begin{cases}
      C_x & \text{ if } \precc x \smallsetminus \upc m = \emptyset, \\
    C_{x'} & \text{ if } \precc x \smallsetminus \upc m=\{x'\},
\end{cases}
\]
    noting that $C$ is a well-defined clopen because we assumed $|\precc x \smallsetminus \upc m | \leq 1$, and that clearly $C$ contains $M_x$.
    It now follows from Lemma \ref{lem compact} that there exists a point $z' \in ]x,z[$ such that $min(\down z' \smallsetminus \down x) \cup M_x \subseteq C.$

    Again using the assumption that $x$ is a down limit, the chain $[x,z']$ is also infinite, and therefore we can apply Lemma \ref{lem inf chains have teeth} twice, obtaining points $(u<v) \in \,]x,z']$ such that $\precc u \smallsetminus \upc m$ and $\precc v \smallsetminus \upc m$ are nonempty sets.
    Since $\X$ is a co-tree, hence its points have at most one immediate successor, the inequality $u < v$ entails that the two previous sets are disjoint.
    It is easy to see that the points $(u<v) \in \up m$ satisfy the following conditions: (i) $col(u)=col(v)$ (since $u,v \in \, ]x,z'] \subseteq [x,z] \subseteq C_x$); (ii) $\precc u \smallsetminus \upc m \neq \emptyset$ (by the definition of $u$); (iii) the set 
    \[
    M \coloneqq \big[(\down v \smallsetminus \down u) \cup (\precc u \smallsetminus \upc m)\big] \cap min(\X)
    \]
    has at least two elements (since it contains the nonempty disjoint sets $\precc u \smallsetminus \upc m$ and $\precc v \smallsetminus \upc m$).
    By Lemma \ref{lem forbidden conditions}, 
    to find our desired contradiction it now suffices to show that $M$ only contains points of the same color.

    Accordingly, we will prove that $M \subseteq C$ (recall that by definition, the set $C$ is equal to a clopen corresponding to a single color).
    Note that the inclusion $min(\down z' \smallsetminus \down x)  \subseteq C$ comes from the definition of the point $z'$, and that we have $(u < v) \in \, ]x,z']$.
    The latter condition clearly entails $(\down v \smallsetminus \down u) \cup (\precc u \smallsetminus \upc m) \subseteq \down z' \smallsetminus \down x$, thus we have
    \[
    M=\big[(\down v \smallsetminus \down u) \cup (\precc u \smallsetminus \upc m)\big] \cap min(\X) \subseteq (\down z' \smallsetminus \down x) \cap min(\X) \subseteq min(\down z' \smallsetminus \down x) \subseteq C. \qedhere
    \]
\end{proof}

\begin{Corollary} \label{corol down are up}
    If $x$ is a down limit, then either $x=m^+$ or $x$ is an up limit. 
    Consequently, we have 
    \[
    \precc x \smallsetminus \upc m = \precc x \subseteq min(\X) \text{ and } 2 \leq | \precc x  | \leq 2^n.\]
\end{Corollary}
\begin{proof}
    Suppose that $x$ is a down limit. 
    By the previous proposition, we have 
    \[
    2 \leq | \precc x \smallsetminus \upc m | \leq 2^n.
    \]
    In particular, $x$ is not minimal, and is thus distinct from $m$.
    Hence $x \in \upc m$.
    Let us now suppose that $x$ is not an up limit.
    By Definition \ref{def up and down limits}, $x \neq m$ yields $\precc x \cap \up m \neq \emptyset$.
    Note that if $\precc x \cap \upc m \neq \emptyset$, Lemma \ref{clop lemma}.(ii) would imply that $\down x $ is clopen.
    But this cannot happen, because $\down x \cap \, ]x,r]= \emptyset$ and $x$ is a down limit, hence also a limit point of the closed chain $[x,r]$.
    Consequently, $\precc x \cap \up m = \precc x \cap \{m\}=\{m\}$, hence $x=m^+$, as desired.

    The second part of the statement is now an immediate consequence of the above display, together with the definitions of up limits and of $m^+$.
\end{proof}

We are finally ready to prove that down limits do not exist.\vs

\noindent \textit{Proof of Theorem} \ref{thm structure of X}.(ii).
    Suppose that $x \in \up m$ is a down limit.
    This entails that $x$ is a limit point of the closed chain $\up x \cap C_x$, where $C_x$ is the clopen corresponding to $col(x)$.
    Hence there must exist $z'\in \upc x$ such that $[x,z'] \subseteq C_x$.
    Moreover, by Corollary \ref{corol down are up} we have
    \[
    \precc x = \precc x \smallsetminus \upc m= \{x_1, \dots , x_t\} \subseteq min(\X),
    \]
    for some $t$ satisfying $2 \leq t \leq 2^n$.
    For each $i \leq t$, we denote by $C_i$ the clopen $C_{x_i}$ corresponding to $col(x_i)$.
    Observe that, by Lemma \ref{lemma order-iso}, distinct points in $\precc x$ must have distinct colors, otherwise this would contradict the Coloring Theorem \ref{coloring thm}.
    Hence the $C_i$'s are pairwise disjoint, and it is clear that they cover $\precc x$.
    Using the notation introduced in Lemma \ref{lem compact}, we have $M_x = \precc x \smallsetminus \upc m = \precc x$ because $x \neq m$.
    Since the $C_i$'s cover $\precc x$, it follows that
    \[
    M_x = \precc x \subseteq C \coloneqq \bigcup_{i=1}^t C_i,
    \]
    and we can apply the aforementioned lemma, yielding a point $z \in \, ]x,z'[\, \subseteq C_x$ satisfying
    \[
    min(\down z \smallsetminus \down x) \cup \precc x = min(\down z \smallsetminus \down x) \cup (\precc x \smallsetminus \upc m) \subseteq C = \bigcup_{i=1}^t C_i.
    \]
    For $i \leq t$, we define $M_i \coloneqq \big(min(\down z \smallsetminus \down x) \cap C_{i} \big) \cup \{x_i\} $, and note that $M_i \subseteq C_i \cap min(\X)$ (recall that $x$ is assumed to be a down limit, hence $min(\down z \smallsetminus \down x)\subseteq \min(\X)$). Set
    \[
    E \coloneqq [x,z]^2 \cup M_1^2 \cup \dots \cup M_t^2 \cup Id_X.
    \]
    
    We will prove that $E$ is a proper bi-E-partition (see Definition \ref{def bi-be}) that only identifies points of the same color, thus contradicting the Coloring Theorem \ref{coloring thm}.
    That $E$ is an equivalence relation is clear, and it is proper because, e.g., we have $xEz$ and $x \neq z$.
    Moreover, since $[x,z] \subseteq C_x$ and $M_i \subseteq C_{i}$ for every $i \leq t$, it is clear that $E$ only identifies points of the same color.
    The proof that $E$ satisfies the up and down conditions is routine and follows closely the argument detailed in the proof of Lemma \ref{lem forbidden conditions} used to establish the same conditions, hence we skip it.
    
    It remains to show that $E$ is refined, that is, for every $u,v \in X$, 
    \[
    \text{ if } (u,v)\notin E, \text{ then there exists an $E$-saturated clopen upset }  U  \text{ s.t. } |U \cap \{u,v\}|=1.
    \]
    Accordingly, let $u,v \in X$ be such that $(u,v)\notin E$, and observe that this is equivalent to demanding that $u \neq v$ and $|[x,z] \cap \{u,v\}| \leq 1 \geq |M_i \cap \{u,v\}|$, for every $i \leq t$.
    We proceed by cases.
    \benbullet
        \item \underline{\textbf{Case:} $u \in M_i$ and $v \in M_j$, for some $i \neq j \leq t$} \\
        We take $U\coloneqq \up C_{i}$, noting that this clopen upset clearly contains $M_i \cup [x,z]$, since $x_i \in M_i \subseteq C_{i}$ and $x_i < x$.
        Now, recall that we established above that $M_k \subseteq C_{k} \cap min(\X)$ and $C_{i} \cap C_{k} = \emptyset$, for all $k \leq t$ such that $k \neq i$.
        Using these conditions, the fact that $U=\up C_i$ is an upset implies that $U$ is disjoint from $M_k$, for all $k \leq t$ such that $k \neq i$.
        Consequently, $U$ is an E-saturated (since it contains $M_i \cup [x,z]$ and is disjoint from all the other $M_k$'s) clopen upset that contains $u \in M_i$ but not $v \in M_j$.\vspace{.3cm}

        \item \underline{\textbf{Case:} $u \in min(\X)$ and $v \in min(\X) \smallsetminus (\bigcup_{i=1}^t  M_i)$} \\
        Recall that $min(\down z \smallsetminus \down x) \cup \precc x \subseteq \bigcup_{i=1}^t C_i$ and $\precc x = \{x_1, \dots , x_t\}$, so it is clear that
        \[
        min(\down z \smallsetminus \down x) \cup \precc x  =\bigcup_{i=1}^t  \big[\big(min(\down z \smallsetminus \down x) \cap C_{i} \big) \cup \{x_i\} \big]= \bigcup_{i=1}^t  M_i,
        \]
        where the second equality holds because of the definitions of the $M_i$'s.
        By the above display, our assumption on $v$ translates to either $v \in  min(\down x) \smallsetminus \precc x$ or $v \in min(\X) \cap (\down r \smallsetminus \down z)$.
        
        If $v \in min(\down x) \smallsetminus \precc x$, then we must have $x \neq m^+$.
        To see this, observe that otherwise the equality $X=\upc m \uplus min(\X)$ established in Proposition \ref{prop diagonal} would entail $min(\down x)=\precc x$, contradicting the assumption that $v \in min(\down x) \smallsetminus \precc x$.
        Notice that this assumption also implies $v <x$, because $v \in \down x$ and since $x$ is not minimal in $\X$ (recall that $\precc x \neq \emptyset$).
        It now follows from Corollary \ref{corol down are up} that $x$ must be an up limit, and as such, from $v <x$ we can certainly find a gap $v < a \prec b <x$ satisfying $(a \prec b ) \in \upc m$.
        By Lemma \ref{clop lemma}.(ii), $\down a$ is a clopen neighbourhood of $v$, which moreover is disjoint from $\bigcup_{i=1}^t  M_i$, by the above display together with $a < x < z$ and $a \not \prec x$.
        It is easy to see that not only does $\up \down a$ contain $\{v\} \cup [x,z]$, but also that the upset $\up \down a$ is disjoint from $\bigcup_{i=1}^t  M_i$, since so is $\down a$ and we know $\bigcup_{i=1}^t  M_i \subseteq min(\X)$.
        Now, since $u$ and $v$ are non-E-equivalent minimal points, they must satisfy $v \nleq u$.
        By the \hyperref[PSA]{PSA}, there exists a clopen upset $V$ separating $v$ from $u$, and since $v < x$ by hypothesis, the fact that $V$ is an upset yields $[x,z] \subseteq V$.
        We take $V' \coloneqq \up \down a \cap V$ as our $E$-saturated clopen upset, as it clearly contains $\{v\} \cup [x,z]$ and is disjoint from $\bigcup_{i=1}^t  M_i \cup \{u\}$.

        On the other hand, if $v \in min(\X) \cap (\down r \smallsetminus \down z)$, then we must have $v\nleq z$.
        By the same argument as above, we know that $v \nleq u$ also holds.
        Thus, by applying the \hyperref[PSA]{PSA} to the two previous inequalities and taking the intersection of the resulting clopen upsets, we can find a clopen upset $V$ which contains $v$ but omits both $z$ and $u$. 
        The clopen $V$ is clearly $E$-saturated, since it is an upset which does not contain $z$, and we have $[x,z] \cup \bigcup_{i=1}^t  M_i \subseteq \down z$. \vspace{.3cm}

        \item \underline{\textbf{Case:} $u \in min(\X)$ and $v \in \upc m$} \\
        Let $v' \coloneqq MIN\{v, x\}$, which exists since both $v$ and $x$ are contained in the chain $\upc m$ (the former by assumption, and the latter because we proved in Corollary \ref{corol down are up} that down limits, such as $x$, must lie in $\upc m$).
        Notice that also $v' \in \upc m$, so our assumption $u \in min(\X)$, together with Corollary \ref{corol diagonal}, yields a clopen upset $V \subseteq \upc m$ which contains $v'$ and omits $u$.
        By the definition of $v'$, we have $\{v\} \cup [x,z]\subseteq V$, while $V \subseteq \upc m = X \smallsetminus min(\X)$ (see Proposition \ref{prop diagonal}) and $\bigcup_{i=1}^t  M_i\subseteq min(\X)$ imply
        \[
        V \cap \big(\bigcup_{i=1}^t  M_i \cup \{u\}\big) = \emptyset.
        \]
        Thus, $V$ satisfies the desired conditions. \vspace{.3cm}

        \item \underline{\textbf{Case:} $u,v \in \upc m$} \\
        Without loss of generality, we assume $u < v$, and proceed by subcases, recalling that $|[x,z] \cap \{u,v\} | \leq 1$ holds by assumption on $u$ and $v$.
        
        If $v<x$, then we take a gap $u \leq a \prec b \leq v$, so by Lemma \ref{clop lemma}.(i), $\up b$ is a clopen upset, which moreover contains $\{v\} \cup [x,z]$ and is disjoint from $\{u\} \cup \bigcup_{i=1}^t  M_i$.

        If $u < x \leq v$, then we take a gap $u \leq a \prec b \leq x$, so by Lemma \ref{clop lemma}.(i), $\up b$ is a clopen upset, which moreover contains $\{v\} \cup [x,z]$ and is disjoint from $\{u\} \cup \bigcup_{i=1}^t  M_i$.

        If $x \leq u \leq z < v$, then we take a gap $z \leq a \prec b \leq v$, so by Lemma \ref{clop lemma}.(i), $\up b$ is a clopen upset, which moreover contains $v$ and is disjoint from $\{u\} \cup [x,z] \cup \bigcup_{i=1}^t  M_i$.

        If $z < u < v$, then we take a gap $u \leq a \prec b \leq v$, so by Lemma \ref{clop lemma}.(i), $\up b$ is a clopen upset, which moreover contains $v$ and is disjoint from $\{u\} \cup [x,z] \cup \bigcup_{i=1}^t  M_i$.

        In all four subcases detailed above (which cover every possibility since in this case, the points $u,v,x,z$ all lie in the chain $\up m$), $V \coloneqq \up b$ is an $E$-saturated clopen upset separating $v$ from $u$, as desired.
    \ebullet

     Using the facts $X = \upc m \uplus min(\X)$ and $\bigcup_{i=1}^t  M_i \subseteq min(\X)$, it is easy to see that the cases detailed above cover all the possibilities for the points $u$ and $v$, since they clearly imply
    \[
    X = \upc m \uplus (\bigcup_{i=1}^t  M_i ) \uplus [min(\X) \smallsetminus \big(\bigcup_{i=1}^t  M_i\big)].
    \]

    We conclude that $E$ is refined, and thus contradicts the Coloring Theorem \ref{coloring thm}.
    It follows that $x$ cannot be a down limit. \qed \vs

Now that we have proved that down limits of $\X$ do not exist, the fact that the diagonal $\up m$ of $\X$ is an infinite well-order follows readily. \vs

\noindent\textit{Proof of Theorem} \ref{thm structure of X}.(iii).
    That the chain $\up m$ is an infinite linear order follows immediately from Corollary \ref{corol diagonal}, so let us suppose that $U$ is a nonempty subset of $\up m$ and show that it must have a least element.
    As $U$ is a nonempty subchain of $\up m$, it has an infimum $u \coloneqq inf(U)$ by Proposition \ref{bi-esa prop}.
    Clearly, $u \in \up m$, and we proved above that $u$ cannot be a down limit, hence either $u=r$ or $u$ has a unique immediate successor $u^+$ in $\up m$.
    If $u=r$, then we must have $U=\{r\}$ and $u=r=MIN(U)$.
    If $u^+$ exists, then $inf(U)=u < u^+$ entails that there exists $v \in U$ such that $v < u^+$.
    As $u,v,u^+ \in \up m$ and $u \prec u^+$, we conclude that $u=v \in U$, thus $u=MIN(U)$, as desired. \qed \vs

\begin{Corollary} \label{corol well-order}
    If $x \in \X \smallsetminus \{r\}$, then $x$ has an unique immediate successor $x^+ \in \upc m$.
\end{Corollary}
\begin{proof}
Notice that, for points in $min(\X) \smallsetminus \{m\}$, this result was already established in Proposition \ref{prop diagonal}.
Furthermore, by the aforementioned proposition, we have $X=\upc m \uplus min(\X)$.
So, to finish this proof, it suffices to prove the statement for points in $\up m \smallsetminus \{r\}$, which follows immediately from Theorem \ref{thm structure of X}.(iii).
\end{proof}

The last two conditions of Theorem \ref{thm structure of X} provide bounds for the number of immediate predecessors of the points in the diagonal $\up m$ which are not up limits.
Equipped with these bounds, we will finally be able to provide a thorough description of the poset structure of $\X$.

\begin{Lemma} \label{lem isolated points}
    If $x \in \X$ is not an up limit, nor an immediate predecessor of an up limit, then $\{x\}$ is a clopen set.
\end{Lemma}
\begin{proof}
    Recall that $X= \upc m \uplus min(\X)$, by Proposition \ref{prop diagonal}.
    We first assume that $x \in \upc m$ is not an up limit, hence there exists $y \in \up m$ such that $y \prec x$, and it follows from Lemma \ref{clop lemma}.(i) that $\up x$ is clopen.
    If $x = r$, then as $\up r$ and $\down r = X$ are both clopen, so is $\up r \cap \down r=\{r\}=\{x\}$.
    If $y\prec x < r$, then combining Corollary \ref{corol well-order} with Lemma \ref{clop lemma}.(iii) implies that $\{x\}$ is clopen.

    Now, suppose that $x \in min(\X)$ is not an immediate predecessor of an up limit.
    By Corollary \ref{corol well-order}, the point $x$ has a unique immediate successor $x\qq +$ contained in $\upc m$.
    By assumption, $x\qq +$ is not an up limit, hence $y \prec x\qq +$ for some $y \in \up m$.
    We established above that $\{x\qq +\}$ must be clopen, hence so is $\downc x\qq +$.
    If $y =m$ then $x\qq += m\qq +$, and the following equalities are easily verified:
    \[
     \downc x^+ = \downc m^+=\precc m^+ = \precc m^+ \smallsetminus \upc m = \precc x\qq + \smallsetminus \upc m.
     \]
    Since $\downc x\qq +$ is clopen, the above display ensures that $\precc x\qq + \smallsetminus \upc m$ is also clopen.
    On the other hand, if $m < y$ then $(y \prec v) \in \upc m$, and it follows from Lemma \ref{clop lemma}.(ii) that $\down y$ must be clopen.
     Thus, the set $\precc x\qq + \smallsetminus \upc m$ is clopen, since so are $\downc x^+$ and $\down y$, and we clearly have $\downc x^+ \smallsetminus \down y= \precc x\qq + \smallsetminus \upc m$.

     In both cases for $y$ detailed above, $\precc x\qq + \smallsetminus \upc m$ is a clopen set of minimal elements (recall that $X= \upc m \uplus min(\X)$), hence it follows from Lemma \ref{lem minimal im pred} that this must be a finite clopen set.
    Therefore, since bi-Esakia spaces are Hausdorff, we can now infer that every point in $\precc x\qq + \smallsetminus \upc m$ is isolated in $\X$, and the result follows.
\end{proof}

\begin{Lemma} \label{lem non up have at most 2 im pred}
    If $x \in \upc m$ is not an up limit, then $|\precc x | \leq 2$.
\end{Lemma}
\begin{proof}
    If $x \in \upc m$ is not an up limit, then there exists $y \in \up m$ such that $y \prec x$.
    Clearly, either $y=m$ or $y \in \upc m$, and both cases yield that $\down y$ is a clopen set; the former by Lemma \ref{lem isolated points}, and the latter by Lemma \ref{clop lemma}.(ii).
    Let us now suppose, with a view to contradiction, that there are distinct points, $x_1$ and $x_2$, contained in $\precc x \smallsetminus \{y\}$.
    Since we assumed that $x$ is not an up limit, both $\{x\}$ and $\{x_1\}$ are clopen by Lemma \ref{lem isolated points}.
    Set $V_a \coloneqq \up x$, $V_b \coloneqq \down y$, $V_c \coloneqq \{x_1\}$, and $V_d \coloneqq X \smallsetminus (V_a \cup V_b \cup V_c)$, and observe that they are nonempty (notice that $V_d$ must contain $x_2$) pairwise disjoint clopen sets (by our previous discussion) that cover $\X$.
    It is not hard to see that the map $f \colon \X \to \F_3$ (see Figure \ref{Fig:the co-trees}) defined by $f[V_e] \coloneqq \{e\}$, for $e \in \{a,b,c,d\}$, is a surjective bi-Esakia morphism.
    Hence $\X \not \models \J(\F_3)$ by the Jankov Lemma \ref{jankov lemma}.
    But this contradicts our assumption that $\X \models LFC$, since $\J(\F_3)$ is an axiom of $LFC$.
    We conclude that $\precc x \smallsetminus \{y\}$ contains at most one point, as desired.
\end{proof}
\vs

\noindent \textit{Proof of Theorem} \ref{thm structure of X}.(iv) \& (v). Notice that, by the definition of up limits \ref{def up and down limits}, condition (iv) is just a particular case of Lemma \ref{lem non up have at most 2 im pred}.
In order to prove (v), we first establish a helpful claim.

\begin{Claim}
    $\precc r \smallsetminus \upc m \neq \emptyset$.
\end{Claim}
\noindent \textit{Proof of the Claim.} If $r$ is an up limit, then this claim follows immediately from Lemma \ref{lem up lims have pred}, so let us suppose that $r$ has an immediate predecessor $y$ in the diagonal of $\X$, noting that since $\up m$ is infinite by Theorem \ref{thm structure of X}.(iii), we must have $y \in \upc m$.
We assume with a view to contradiction that $\precc r \smallsetminus \upc m = \emptyset$.
As $y \in \upc m$, it follows from Corollary \ref{corol diagonal} that there exists a clopen upset $U$ satisfying $y \in U \subseteq \upc m$.
Observe that: $V_a \coloneqq \{r\}$ is a nonempty clopen set by Lemma \ref{lem isolated points}; $V_b \coloneqq U \smallsetminus \{r\}$ is a nonempty clopen set since $U$ is a clopen neighbourhood of $y$, a point that is distinct from $r$; and that $V_c \coloneqq \down V_b \smallsetminus V_b$ is a clopen set, since so are $\down V_b$ and $V_b$, which is moreover nonempty, because $m \notin V_b \subseteq U \subseteq \upc m$ and $y \in \upc m \cap  V_b$ entails $m \in \down V_b$.
Furthermore, the sets $V_a$, $V_b$, and $V_c$ are pairwise disjoint clopens that cover $\X$ and satisfy the following relations:
\[
V_a=\up V_b \smallsetminus V_b \text{ and } \down V_a \smallsetminus V_a = V_b \uplus V_c \text{ and } \up V_c \smallsetminus V_c = V_a \uplus V_b
\]
(this is very easily seen by picturing $\X$ and the three clopens, making use of our assumption that $\precc r \smallsetminus \upc m = \emptyset$).
It is routine to check that the map $f \colon \X \to \F_1$ (see Figure \ref{Fig:the co-trees}) defined by $f[V_e]\coloneqq \{e\}$, for $e\in \{a,b,c\}$, is a surjective bi-Esakia morphism, so the Jankov Lemma \ref{jankov lemma} yields $\X \not \models \J(\F_1)$.
But this contradicts our assumption that $\X \models LFC$, since $\J(\F_1)$ is an axiom of $LFC$. \qed \vs

Let us now prove Theorem \ref{thm structure of X}.(v).
Suppose that $x\in \up m \smallsetminus \{m,m^+\}=\upc m^+$ is not an up limit. We need to establish
\[
|\precc x \cap \upc m|=1=|\precc x \smallsetminus \upc m|.
\]
By the definition of up limits \ref{def up and down limits}, our assumptions on $x$ imply that there exists $y \in \upc m$ such that $y \prec x$.
As $\upc m$ is a chain, it is clear that $\precc x \cap \upc m = \{y\}$, hence it remains to show that $|\precc x \smallsetminus \upc m| =1$. 
But $x$ falls under the conditions of Lemma \ref{lem non up have at most 2 im pred}, so we know that $|\precc x| \leq 2$.
Thus, to establish the desired equality, it suffices to prove that $\precc x \smallsetminus \upc m$ is nonempty, because $\precc x \cap \upc m = \{y\}$.

We assume $\precc x \smallsetminus \upc m = \emptyset$ and arrive at a contradiction.
By Corollary \ref{corol well-order}, either $x$ is the co-root $r$ of $\X$ or $x^+\in \upc m$ exists.
Since the case $x=r$ follows immediately from the previously established claim, we suppose that $ x^+\in \upc m$ exists, noting that we now have $(y \prec x \prec x^+)\in \upc m$, which entails that $\{x\}$ is clopen by Lemma \ref{clop lemma}.(iii).
By the aforementioned claim and by Theorem \ref{thm structure of X}.(iii), the set $H\coloneqq \{ v \in \up x\qq + \colon \precc v \smallsetminus \upc m \neq \emptyset \} \subseteq \up m$ is a nonempty subset of a well-order, and thus has a least element $z$.
Notice that we have $x\qq + \leq z$, hence $x < z$ follows.

We will show that the nonempty interval $[x,z[$ is a finite clopen set.
For suppose that it is infinite. 
Then there must exist $v \in [x,z[$ satisfying $|[x,v]| \geq 2\qq n$, so Lemma \ref{lem inf chains have teeth} contradicts the definition of $z$ as the least element of $H$.
Therefore, the chain $[x,z[$ must be finite and we infer that $z$ has a unique immediate predecessor contained in $[x,z[\, \subseteq \upc m$, which we denote by $z'$.
Since $(z' \prec z) \in \upc m$, Lemma \ref{clop lemma} entails that both $\down z'$ and $\up z$ are clopen sets.
As $\{x\}$ is also clopen by above, so is $\up x$, and we conclude that $[x,z[ \, =[x,z']=\up x \cap \down z'$ is a finite clopen set, as desired.

Set $V_a \coloneqq \up z$, $V_{a'}\coloneqq X \smallsetminus(\down z' \cup \up z)$, and $V_b \coloneqq [x,z']$, noting that our previous discussion ensures that they are pairwise disjoint clopen sets, and that both $V_a$ and $V_b$ are nonempty.
Moreover, using $z \in H$ and the fact that $z'$ is the sole element of $\precc z \cap \upc m$, it is easy to see that both $V_{a'}\neq \emptyset$ and $V_a = \up V_{a'}\smallsetminus V_{a'}$ hold true.
Now, recall that $(y \prec x \prec x^+)\in \upc m$, so in particular we have $y \in \upc m$.
It follows from Lemma \ref{corol diagonal} that there is a clopen upset $U$ such that $y \in U \subseteq \upc m$.
Hence the fact that $\up x$ is clopen (see above) implies that $V_c \coloneqq U \smallsetminus \up x$ is clopen as well, and clearly so is $V_d \coloneqq \down V_c \smallsetminus V_c$.
Moreover, these two clopens are nonempty: the former contains $y$ because $y \in U$ and $y < x$, and the latter contains $m$ because $m < y \in V_c \subseteq U$ and $m \not \in U \subseteq \upc m$.
    It is not hard to see that $\{V_a, V_{a'}, V_b, V_c ,V_d\}$ is a family of nonempty pairwise disjoint clopen sets that moreover covers $\X$ (this is easily seen by picturing $\X$ and the five clopens, while making crucial use of the fact that the definition of $z$ ensures that $\precc v \smallsetminus \upc m = \emptyset$, for every $v \in [x,z']$).
    It is routine to check that the map $f \colon \X \to \F_2$ (see Figure \ref{Fig:the co-trees}) defined by $f[V_e] \coloneqq \{e\}$, for $e \in \{a,a',b,c,d\}$, is a surjective bi-Esakia morphism, so $\X \not \models \J(\F_2)$ follows from the Jankov Lemma \ref{jankov lemma}.
    As $\J(\F_2)$ is an axiom of $LFC$ and we assumed that $\X$ validates this logic, we found the desired contradiction. \qed \vs

We will finally describe the poset structure of $\X$.
In Figure \ref{Fig:general structure}, points are painted blue when they might be identified with other points.
More precisely, $y_0$ is painted blue because we might have that $x_0=y_0$, depending on the particular bi-Esakia space $\X$. 
It can also happen that $n'=0$, in which case we have $\eta' = \eta$ and $ t_{\eta'}=1$ and $y_\eta = y_{\eta'}^1$.

\begin{Theorem} \label{Thm poset structure}
    Let $n$ be a positive integer and $\A$ an infinite $n$-generated SI bi-Heyting algebra.
    If $\A$ validates
    \[
    LFC = \lc + \beta(\F_0) +\J(\F_1) + \J(\F_2) +\J(\F_3),
    \]
    then the underlying poset of $\A_*$ looks like the poset represented in Figure \ref{Fig:general structure}, where: $\eta$ is an infinite ordinal such that $|\eta| \leq 2 \qq {\aleph_0}$; the ordinals $n'$ and $\eta'$ satisfy $n'< \omega \leq \eta' \leq \eta = \eta' + n'$;
    and $2\leq t_\lambda \leq 2\qq n$, for every limit ordinal $\lambda$ such that $\omega \leq \lambda \leq \eta'$.
\end{Theorem}

    \begin{figure}[h]
\centering
\begin{tabular}{c}
\begin{tikzpicture}
    \tikzstyle{point} = [shape=circle, thick, draw=black, fill=black , scale=0.35]
    \tikzstyle{bpoint} = [shape=circle, thick, draw=blue, fill=blue , scale=0.35]
    \node [label=left:{$x_0$}] (0) at (0,0) [point] {};
    \node [label=right:{$y_0$}] (0') at (1,-1) [bpoint] {};
    \node [label=right:{$y_1$}] (1') at (1.75,-.25) [point] {};
    \node [label=left:{$x_1$}] (1) at (0.75,0.75) [point] {};
    \node [label=left:{$x_\omega$}] (omega) at (2,2) [point] {};
    \node [label=right:{$y_\omega^1$}] (a) at (2.7,0.7) [point] {};
    \node [label=right:{$y_\omega^{t_\omega}$}] (b) at (3.3,1.3) [point] {};
    \node [label=left:{$x_{\omega +1}$}] (omega+1) at (3,3) [point] {};
    \node [label=right:{$y_{\omega + 1}$}] (omega+1') at (4,2) [point] {};
    \node [label=left:{$x_{\lambda}$}] (alpha) at (4.25,4.25) [point] {};
    \node [label=right:{$y_\lambda^1$}] (a2) at (4.95,2.95) [point] {};
    \node [label=right:{$y_\lambda^{t_\lambda}$}] (b2) at (5.55,3.55) [point] {};
    \node [label=left:{$x_{\lambda+1}$}] (alpha1) at (5.25,5.25) [point] {};
    \node [label=right:{$y_{\lambda+1}$}] (alpha1') at (6.25,4.25) [point] {};
    \node [label=left:{$x_{\eta'}$}] (eta') at (6.5,6.5) [point] {};
    \node [label=right:{$y_{\eta'}^1$}] (a3) at (7.2,5.2) [point] {};
    \node [label=right:{$y_{\eta'}^{t_{\eta'}}$}] (b3) at (7.8,5.8) [point] {};
    
    \node [label=left:{$x_{\eta'+1}$}] (eta'+1) at (7.5,7.5) [bpoint] {};
    \node [label=right:{$y_{\eta'+1}$}] (eta'+1') at (8.5,6.5) [bpoint] {};
    \node [label=above:{$x_{\eta'+n'}=x_\eta$}] (eta) at (8.5,8.5) [bpoint] {};
    \node [label=right:{$y_{\eta'+n'}=y_\eta$}] (eta1) at (9.5,7.5) [bpoint] {};

    \color{blue}
    \draw  (0')--(0);  \color{black}
    \draw (0)--(1)--(1');
    \draw[dotted] (1)--(omega)--(alpha);
    \draw (omega)--(omega+1)--(omega+1');
    \draw (a)--(omega)--(b);
    \draw[dotted] (a)--(b);
    \draw[dotted] (a2)--(b2);
    \draw (a2)--(alpha)--(b2);
    \draw (alpha)--(alpha1)--(alpha1');
    \draw (a3)--(eta')--(b3);
    \draw[dotted] (alpha1)--(eta);
    \draw[dotted] (a3)--(b3); \color{blue}
    \draw (eta')--(eta'+1)--(eta'+1');
    \draw (eta)--(eta1);
\end{tikzpicture}
\end{tabular}
\caption{The poset $\A_*$.}
\label{Fig:general structure}
\end{figure}

\begin{proof}
    Let $\A$ satisfy the conditions of the statement.
    Without loss of generality, we assume that $\A_*=\X$, where $\X$ is the co-tree we have been working with throughout this section.
    We start by noting that the language of $\bipc$ is countable, so the assumption that $\A$ is finitely generated implies that $\A$ is also countable.
    It follows from the bi-Esakia duality that $|X| \leq 2^{|A|} =2^{\aleph_0}$.
    As the diagonal $\up m$ of $\X$ is an infinite well-order by Theorem \ref{thm structure of X}.(iii), we have that $\up m$ is order isomorphic to an infinite ordinal $\eta$, hence the previous inequality entails that $|\eta| \leq 2 ^{\aleph_0}$.

    Let $f\colon \eta \to \up m$ be the order isomorphism coming from Theorem \ref{thm structure of X}.(iii). 
    If $\precc m^+ = \{m\}$, we set $y_0\coloneqq m$ and $x_0 \coloneqq m^+$.
    If there exists $m_1 \in \precc m^+ \smallsetminus \{m\}$, then it follows from Theorem \ref{thm structure of X}.(iv) that $m_1$ is the unique immediate predecessor of $m^+$ that is distinct from $m$.
    In this case, we define both $y_0$ and $x_0$ to be the point $m$, and set $y_1 \coloneqq m_1$ and $x_1 \coloneqq m^+$.
    
    Let $x \in \upc m^+$.
    If $x$ is not an up limit, then it has a unique immediate predecessor $z$ in $\upc m$.
    Suppose we have already defined $x_\alpha \coloneqq z$, for some ordinal $\alpha \in ]0,\eta[\,$.
    We set $x_{\alpha +1}\coloneqq x$, and let $y_{\alpha+1}$ be the unique point contained in the set $\precc x \smallsetminus \upc m$, which exists by Theorem \ref{thm structure of X}.(v).
    On the other hand, if $x$ is an up limit, then Theorem \ref{thm structure of X}.(i) implies that $\precc x = \{z_1, \dots ,z_t\}$, for some $2 \leq t \leq 2^n$.
    Moreover, there exists a unique ordinal $\lambda \leq \eta$ such that $f(\lambda)=x$.
    We set $x_\lambda \coloneqq x$ and $y_\lambda^i \coloneqq z_i$, for $i \leq t$.
    Since Proposition \ref{prop diagonal} ensures that $X=\upc m \uplus min(\X)$ and that every point in $min(\X)\smallsetminus\{m\}$ has an immediate successor in $\upc m$, the above argument is sufficient to prove the statement.
\end{proof}

\subsection{The Finite Model Property}

We will continue working with a fixed but arbitrary infinite bi-Esakia co-tree $\X$ which validates $LFC$, and whose dual is $n$-generated, for some $n\in \ZZ^+$.
By Theorem \ref{Thm poset structure}, we know the structure of the underlying poset of $\X$.
In particular, we are working with fixed but arbitrary ordinals $n'$, $\eta'$ and $\eta$ satisfying $n'< \omega \leq \eta' \leq \eta = \eta' + n'$.
In order to simplify some statements, we shall use the convention that $0$ is not a limit ordinal.

If $z$ is a point in $\X$, we define the \textit{ordinal of $z$} as 
\[
ord(z) \coloneqq 
\begin{cases}
    \alpha & \quad \text{ if } z \in \{x_\alpha, y_\alpha\} \text{ and } \alpha \text{ is a nonlimit ordinal},\\
    \lambda & \quad \text{ if } z \in \{x_\lambda, y_\lambda^1, \dots , y_\lambda^{t_\lambda}\} \text{ and $\lambda$ is a limit ordinal}.
\end{cases}
\]
By the structure of $\X$, it is clear that the map $ord$ is order preserving on $X$, and order invariant on $\up x_0$ (i.e., $z \leq x$ iff $ord(z) \leq ord(x)$, for every $z,x \in \up x_0$). In other words, when restricted to $\up x_0$, the map $ord$ is an order embedding.
These facts will be used without further reference in what follows. \vs

The following lemma translates some results of the previous section into our new ``ordinal terminology".

\begin{Lemma} \label{lem limit and nonlimit}
    Let $\alpha, \lambda \leq \eta$ be ordinals.
    If $\alpha$ is a nonlimit ordinal, then $\{x_\alpha\}$ and $\{y_\alpha\}$ are both clopen sets.
    If $\lambda $ is a limit ordinal, then $\up x_\lambda $ is not a clopen set, and $x_\lambda$ is a limit point of the closed chain $[x_\beta,x_\lambda]$, for every ordinal $\beta < \lambda$.
\end{Lemma}
\begin{proof}
    Let $\alpha, \lambda \leq \eta$ be ordinals.
    If $\alpha$ is a nonlimit ordinal, then by the definition of the map $ord$ and by structure of $\X$ (see Theorem \ref{Thm poset structure}), it is clear that $x_\alpha$ cannot be an up limit (see Definition \ref{def up and down limits}). 
    It now follows from Lemma \ref{lem isolated points} that both $\{x_\alpha\}$ and $\{y_\alpha\}$ are clopen.

    On the other hand, if $\lambda$ is a limit ordinal then $x_\lambda$ must be an up limit.
    Consequently, it follows from Remark \ref{rem top properties} that $x_\lambda$ is a limit point of the closed chain $[z,x_\lambda]$, for every $z\in [y_0,x_\lambda[\,$.
    Since the map $ord$ is order invariant on $\up x_0$, we have that $x_\lambda$ is a limit point of the closed chain $[x_\beta,x_\lambda]$, for every $\beta < \lambda$.
    Finally, notice that, if $\up x_\lambda$ was clopen, then taking $\beta < \lambda$ (which exists, since $\lambda$ is a limit ordinal) would yield $[x_\beta,x_\lambda[\, \cap \up x_\lambda=\emptyset$, thus contradicting the fact that $x_\lambda$ is a limit point of $[x_\beta,x_\lambda]$.
\end{proof}

\begin{Lemma} \label{lem nonempty clopens}
    If $U$ is a nonempty clopen subset of $\X$, then $z \coloneqq MIN(\up y_0 \cap \up U)$ exists, $ord(z)$ is a nonlimit ordinal, and $\{x_\alpha, y_\alpha \} \cap U \neq \emptyset$.
\end{Lemma}
\begin{proof}
    Let $U$ be a nonempty clopen subset of $\X$, and recall that $\{y_0\}$ is a clopen set, by Lemma \ref{lem limit and nonlimit}.
    By the structure of $\X$ (see Theorem \ref{Thm poset structure}), the set $ \up y_0 \cap\up U $ is a nonempty clopen chain in a bi-Esakia space, and therefore has a least element $z$ by Proposition \ref{bi-esa prop}. Clearly, $\up z =\up y_0 \cap \up U $, so it follows that $\up z$ is clopen.
    By Lemma \ref{lem limit and nonlimit}, $\alpha \coloneqq ord(z) $ cannot be a limit ordinal.
    Moreover, by the structure of $\X$, it is now easy to see that the definition of $z$ implies $x_\alpha \in U$ or $y_\alpha \in U$.
\end{proof}

\begin{Theorem} \label{thm finite submodel}
    If $\X$ refutes a formula $\varphi$, then $\X$ has a finite subposet $\Y$ that also refutes $\varphi$ and is either isomorphic to a finite comb, or to a finite hcomb.
\end{Theorem}
\begin{proof}
    Let $V\colon \mathbf{Fm} \to ClopUp(\X)$ be a valuation on $\X$, and suppose that $\M \coloneqq (\X,V) \not \models \varphi$.
    It follows that $X \smallsetminus V(\varphi)$ is a nonempty clopen set.
    By Lemma \ref{lem nonempty clopens}, there exists $u \in X \smallsetminus V(\varphi)$ such that $ord(u) = \alpha$, for some nonlimit ordinal $\alpha \leq \eta$.
    We define:
    \benbullet
        \item $S(\varphi)\coloneqq \{ \psi \in Fm \colon \psi \text{ is a subformula of } \varphi \}$; \vspace{.2cm}

        \item $S_{\gets}(\varphi) \coloneqq \{ \theta \in S(\varphi) \colon  \text{$\gets$ is the main connective of } \theta \};$ \vspace{.2cm}

        \item $S_{\to}(\varphi) \coloneqq \{\sigma \in S(\varphi) \colon \text{$\to$ is the main connective of } \sigma \};$ \vspace{.2cm}

        \item for $\theta=\phi \gets \psi \in S_{\gets}(\varphi)$,
        \[MinOrd(\theta) \coloneqq \begin{cases}
        ord\big( MIN(\up y_0 \cap V(\theta)) \big) & \quad \text{if } V(\theta) \neq \emptyset, \\
        \alpha & \quad \text{otherwise;}
        \end{cases}\]
        

        \item for $\sigma = \phi \to \psi \in S_{\to}(\varphi)$,
        \[MaxOrd(\sigma) \coloneqq \begin{cases}
        ord \big( MAX ( \up y_0 \cap (V(\phi) \smallsetminus V(\psi)))\big)  & \quad \text{if } \up y_0 \cap \big( V(\phi) \smallsetminus V(\psi) \big) \neq \emptyset, \\
        \eta^2 & \quad \text{otherwise\footnotemark.}
        \end{cases}\]
        \footnotetext{The use of $\eta^2$ in this definition is somewhat arbitrary: any ordinal strictly larger than $\eta$ would suffice for our goal.}
    \ebullet
    
    We claim that $MinOrd$ is a well-defined map, whose images are always nonlimit ordinals. 
    To see this, let $\theta = \phi \gets \psi \in S_{\gets}(\varphi)$ and notice that the case where $V(\theta) = \emptyset$ is clear, by the definition of $\alpha$.
    On the other hand, if $V(\theta)=\up (V(\phi) \smallsetminus V(\psi)) \neq \emptyset$, then also the clopen set $V(\phi) \smallsetminus V(\psi)$ must be nonempty, and the claim follows by Lemma \ref{lem nonempty clopens}. 
    Furthermore, note that in this case the following equality follows easily from the structure of $\X$:
    \begin{equation} \label{minord}
        ord\big( MIN(\up y_0 \cap V(\theta)) \big) = MIN\{\delta \in Ord \colon \{x_\delta,y_\delta\}\cap \big(V(\phi) \smallsetminus V(\psi)\big) \neq \emptyset\}.
    \end{equation}
    This is because if $\beta = ord\big( MIN(\up y_0 \cap V(\theta)) \big)$, then $x_\beta =MIN(\up x_0 \cap V(\theta))$.
    As $V(\theta)=\up(V(\phi) \smallsetminus V(\psi))$, it now follows that $\{x_\beta, y_\beta\} \cap V(\phi) \smallsetminus V(\psi) \neq \emptyset$, and also that $\{x_\delta, y_\delta\} \cap V(\phi) \smallsetminus V(\psi) = \emptyset$, for every $\delta < \beta$.
    
    On the other hand, if $\beta = MIN\{\delta \in Ord \colon \{x_\delta,y_\delta\}\cap \big(V(\phi) \smallsetminus V(\psi)\big) \neq \emptyset\}$, then either $y_0= MIN(\up y_0 \cap V(\theta))$ and we have $\beta = 0 = ord\big(MIN(\up y_0 \cap V(\theta))\big)$, or $y_0 < MIN(\up y_0 \cap V(\theta))$, in which case we have $0 < \beta$ and clearly $x_\beta = MIN(\up y_0 \cap V(\theta)) $, hence $\beta = ord\big( MIN(\up y_0 \cap V(\theta)) \big)$.
    This proves (\ref{minord}).

    We also have that $MaxOrd$ is a well-defined map, since if $\up y_0 \cap \big( V(\phi) \smallsetminus V(\psi) \big)$ is nonempty, then it is a nonempty closed chain in $\X$, and is therefore bounded by Proposition \ref{bi-esa prop}. 
    Thus, $MaxOrd(\phi \to \psi)=ord(z)$, where $z$ is the greatest element of the aforementioned closed chain.

    We will now begin our construction of the subposet $\Y$. Set 
    \[
    S_{\to}\qq {lim}(\varphi) \coloneqq \{ \sigma \in S_{\to}(\varphi) \colon MaxOrd(\sigma) \text{ is a limit ordinal}\},
    \]
    and 
    \begin{align*}
        Y_0 \coloneqq  \{x_\alpha,y_\alpha\} & \cup \{x_\delta, y_\delta \colon \exists \, \theta \in S_{\gets}(\varphi) \, \big( MinOrd(\theta) = \delta\big) \} \\
        & \cup \{x_\delta, y_\delta \colon \exists \, \sigma \in S_{\to}(\varphi) \smallsetminus S_{\to}\qq {lim}(\varphi) \, \big(MaxOrd(\sigma)= \delta \leq \eta\big)\},
    \end{align*}
    noting that $Y_0$ must be finite because so is $S(\varphi)$.
    By the same reason, the set $S_{\to}\qq {lim}(\varphi)$ is also finite, and we can thus find a finite enumeration of
    \[
    \Lambda \coloneqq \{\lambda \leq \eta \colon \exists\, \sigma \in S_{\to}\qq {lim}(\varphi) \, \big( MaxOrd(\sigma)= \lambda \big)\}=\{\lambda_1, \dots , \lambda_k\},
    \]
    for some $k \in \omega$.
    We assume, without loss of generality, that $\lambda_1 < \dots < \lambda_k$.
    For each $i \leq k$, it follows from the definitions of $S_{\to}\qq {lim}(\varphi)$, of $\Lambda$, and of $MaxOrd$ that the set
    \[
    U_i\coloneqq \bigcap \{V(\phi)\smallsetminus V(\psi) \colon \phi \to \psi \in S_{\to}\qq {lim}(\varphi) \text{ and } MaxOrd(\phi \to \psi)=\lambda_i\}
    \]
    is a finite intersection of clopens, all of which have $x_{\lambda_i}$ as their unique maximal element that lies in $\up y_0$.
    Hence, $U_i$ is also a clopen set that satisfies $x_{\lambda_i}=MAX( \up y_0 \cap U_i)$. 
    
    Since $Y_0$ is finite by construction, the following ordinals must exist:
    \[
        \beta_1 \coloneqq  MAX\big( \{0\} \cup \{\delta \in Ord \colon x_\delta \in Y_0 \text{ and } \delta < \lambda_1\} \big);
    \]
    and, for $1 < i \leq k$,
    \[
        \beta_i \coloneqq  MAX \big( \{\lambda_{i-1} \} \cup  \{\delta \in Ord \colon x_\delta \in Y_0 \text{ and } \delta < \lambda_i\} \big).
    \]
    Let $i \leq k$, and observe that since $\beta_i+1$ is nonlimit ordinal, $\up x_{\beta_i+1}$ is a clopen upset by Lemma \ref{lem limit and nonlimit}.
    Since $\lambda_i$ is a limit ordinal greater than $\beta_i$ by definition, we have $\beta_i +1 < \lambda_i$, so the upset $\up x_{\beta_i+1}$ contains $x_{\lambda_i}$, and by our comment above, so does $U_i$.
    Therefore, $\up x_{\beta_i+1} \cap U_i$ is a nonempty clopen set, and it follows from Lemma \ref{lem nonempty clopens} that there exists $z \coloneqq MIN\big(\up y_0 \cap \up (\up x_{\beta_i+1} \cap U_i)\big)$ such that $\gamma_i \coloneqq ord(z)$ is a nonlimit ordinal and $\{x_{\gamma_i},y_{\gamma_i}\} \cap (\up x_{\beta_i+1} \cap U_i)\neq \emptyset$.
    By the structure of $\X$ (or, in particular, by the structure of $\up x_{\beta_i+1}$), we know that $y_{\gamma_i} \not \in \up x_{\beta_i+1}$, hence we must have $x_{\gamma_i}\in \up x_{\beta_i+1} \cap U_i \subseteq \up y_0 \cap U_i$.
    As $x_{\lambda_i}=MAX( \up y_0 \cap U_i)$ by above, this proves $x_{\gamma_i} \leq x_{\lambda_i}$.
    But since $\lambda_i$ is a limit ordinal, while $\gamma_i$ is not, we conclude
    \[
    x_{\beta_i} < x_{\beta_i+1} \leq x_{\gamma_i} < x_{\lambda_i}.
    \]
    Because of this, and using the definition of the $\beta_i$'s, the following inequality is clear:
    \[
    \gamma_1 < \lambda_1 \leq \beta_2 < \beta_2 + 1 \leq \gamma_2 < \lambda_2 < \dots \leq \beta_{k-1} < \beta_{k-1} + 1 \leq \gamma_k < \lambda_k.
    \]
    In particular, we have 
    \begin{equation}\label{ineq for claim}
        \gamma_1 < \lambda_1 < \gamma_2 < \lambda_2 < \dots < \gamma_k < \lambda_k.
    \end{equation}
    
    Finally, we set $\Y \coloneqq (Y, \tau, R)$, where $Y \coloneqq Y_0 \cup \{x_{\gamma_i},y_{\gamma_i} \colon 1 \leq i \leq k\}$, $\tau$ is the discrete topology, and $R\coloneqq \,\leq \cap \, Y^2$.
    By construction, $\Y$ contains finitely many points, all of which have nonlimit ordinals.
    By the structure of $\X$ (see Theorem \ref{Thm poset structure}), it is not hard to see that the bi-Esakia co-tree $\Y$ is isomorphic to the $(\frac{|Y|}{2})$-comb if $|Y|$ is even, or isomorphic to the $(\frac{|Y|-1}{2})$-hcomb if $|Y|$ is odd (the latter case only happens when $y_0$ is identified with $x_0$ in $\X$, and $x_0\in Y$).
    
    It remains to show that $\Y$ refutes $\varphi$.
    To this end, let $V{\restriction}\colon \mathbf{Fm}\to ClopUp(\Y)=Up(\Y)$ be the unique valuation satisfying $V{\restriction}(p)\coloneqq V(p) \cap Y$, for every propositional variable $p$, and consider the bi-Esakia model $\mathfrak{N}\coloneqq (\Y, V{\restriction})$.

    \begin{Claim}
         If $z \in \mathfrak{N}$ and $\theta \in S(\varphi)$, then 
    \[
    \M, z \models \theta \iff \mathfrak{N}, z \models \theta.
    \]
    \end{Claim}
    
    \noindent\textit{Proof of the Claim.} We will prove this by induction on the implicative degree (see Section \ref{sec bi-ipc}) of the subformulas of $\varphi$.
    If $ipd(\theta)=0$, the result is immediate by the definitions of implicative degree and of $V{\restriction}$.
    Now, suppose that $\theta \in S(\varphi)$ is such that $ipd(\theta)=i+1$, and also that the statement holds true for every subformula of $\varphi$ whose implicative degree is not greater than $i$.

    We start by considering the cases where the main connective of $\theta$ is an implication or co-implication, i.e., when $\theta = \phi \gets \psi$ or $\theta = \phi \to \psi$, for some $\phi, \psi \in S(\varphi)$.
    By the definition of $ipd$, we have $ipd(\phi),ipd(\psi) \leq i$, hence both formulas fall under our induction hypothesis. In other words, we have 
    \[
    V{\restriction}(\phi)=V(\phi) \cap Y \text{ and } V{\restriction}(\psi)=V(\psi) \cap Y.
    \]
    We will show that $\M, z \models \theta \iff \mathfrak{N}, z \models \theta$, for every $z \in \mathfrak{N}$.


    \benbullet
        \item \underline{\textbf{Case:} $\theta= \phi \gets \psi$}
        
        If $\M, z \not \models \phi \gets \psi$, i.e., if $
        \forall y \leq z\, \big( \M, y \not \models \phi \text{ or } \M, y \models \psi \big),$
        then our induction hypothesis and the definition of $R= \, \leq \cap \, Y^2$ clearly entail $\mathfrak{N}, z \not \models \phi \gets \psi$.

        Suppose now $\M , z \models \phi \gets \psi$, i.e., that $z \in V(\phi \gets \psi)=\up(V(\phi) \smallsetminus V(\psi))$.
        The case $z \in V(\phi) \smallsetminus V(\psi)$ is immediate from our induction hypothesis, so let us assume that $z \notin V(\phi) \smallsetminus V(\psi)$.
        Since $z \in \up (V(\phi) \smallsetminus V(\psi))$, $z$ cannot be a minimal point of $\M$, hence $z \in \upc y_0$ by the structure of $\X$ (see Theorem \ref{Thm poset structure}).
        Moreover, by the definitions of $ord$ and $MinOrd$, it is now clear that, for $\delta \coloneqq MinOrd(\phi \gets \psi)$, the following holds:
        \[
        \delta = ord(MIN(\up y_0 \cap V(\phi \gets \psi) ) \leq ord(z) \text{ and } x_\delta \leq z.
        \]
        By construction, we have $\{x_\delta, y_\delta\} \subseteq Y_0 \subseteq Y$.
        So our assumption that $z\in Y$, together with the fact $y_\delta \leq x_\delta \leq z$, now entails $y_\delta R x_\delta R z$. 
        Furthermore, by the equality (\ref{minord}), we know 
        \[
        \{x_\delta, y_\delta\} \cap (V(\phi)\smallsetminus V(\psi)) \neq \emptyset,
        \]
        hence our induction hypothesis yields $\{x_\delta, y_\delta\} \cap (V{\restriction}(\phi)\smallsetminus V{\restriction}(\psi)) \neq \emptyset$.
        Using $y_\delta R x_\delta R z$, we can now conclude that $z \in R[V{\restriction}(\phi)\smallsetminus V{\restriction}(\psi)]$, i.e., that $\mathfrak{N}, z \models \phi \gets \psi,$ as desired. \vs

        \item \underline{\textbf{Case:} $\theta= \phi \to \psi$}
        
        Firstly, we note that if $\M, z \models \phi \to \psi$, then the definition of $R$ and our induction hypothesis already imply that $\mathfrak{N}, z \models \phi \to \psi$.

        Next, we assume $\M, z \not \models \phi \to \psi$, i.e., that $z \in \down (V(\phi)\smallsetminus V(\psi))$.
        The case $z \in V(\phi)\smallsetminus V(\psi)$ is an easy consequence of the induction hypothesis, so let us assume otherwise, i.e., that $z \notin V(\phi)\smallsetminus V(\psi)$ but there exists a point $w\in \upc z$ such that $w \in V(\phi)\smallsetminus V(\psi)$.
        As $w$ is not a minimal point of $\X$, it follows from the poset structure of this bi-Esakia space (see Theorem \ref{Thm poset structure}) that $w \in \upc y_0$.
        We now have that $w$ is contained in $\up y_0 \cap (V(\phi)\smallsetminus V(\psi))$, hence this set is nonempty.
        By the definition of $MaxOrd$, it follows
        \[
        \delta \coloneqq MaxOrd(\phi \to \psi) = ord\big(MAX(\up y_0 \cap (V(\phi)\smallsetminus V(\psi))\big) \leq \eta.
        \] 
        Furthermore, as $w \in \upc y_0$ and $w \in V(\phi)\smallsetminus V(\psi)$, the definition of $\delta$ forces $x_\delta \in V(\phi)\smallsetminus V(\psi)$.
        Thus, the inequality $z < w \leq x_\delta$ holds in $\X$.

        If $\delta$ is a nonlimit ordinal, then $\delta \leq \eta$ entails $x_\delta \in Y_0 \subseteq Y$ by construction.
        As we are working with the assumption that $z \in Y$, the previous inequality now yields $z R x_\delta$.
        Since, by induction hypothesis, we also have $x_\delta \in Y\cap (V(\phi)\smallsetminus V(\psi))$ i.e., that $x_\delta \in V{\restriction}(\phi) \smallsetminus V{\restriction} (\psi)$, we conclude $\mathfrak{N}, z \not \models \phi \to \psi$.

        We now suppose that $\delta$ is a limit ordinal.
        By the definitions of $MaxOrd$ and of $\delta$, it follows that $\phi \to \psi\in S_{\to}\qq {lim}(\varphi) = \{ \sigma \in S_{\to}(\varphi) \colon MaxOrd(\sigma) \text{ is a limit ordinal}\}$ and $\delta \leq \eta$.
        Recall our definition of the set 
        \[
        \Lambda = \{\lambda \leq \eta \colon \exists\, \sigma \in S_{\to}\qq {lim}(\varphi) \, \big( MaxOrd(\sigma)= \lambda \big)\}=\{\lambda_1, \dots , \lambda_k\},
        \]
        and observe that there must exist an $i\leq k$ such that $\delta = \lambda_i$.
        Thus, the inequality $z < w \leq x_\delta$ now implies that $z < x_{\lambda_i}$ and $ord(z) \leq \lambda_i$.
        In fact, since $\lambda_i$ is a limit ordinal and, by assumption, $z$ is contained in $Y$, a set whose points all have nonlimit ordinals, we can infer $ord(z) < \lambda_i$.
        
        We will show that $z \leq x_{\gamma_i}$.
        For suppose $z \in Y_0$.
        Then $x_{ord(z)}\in Y_0$ by construction, and since $ord(z) < \lambda_i$ by above, it follows from the definitions of $\beta_i$ and of $\gamma_i$ that $ord(z) \leq \beta_i < \gamma_i$.
        Thus, in this case we have $z < x_{\gamma_i}$.
        On the other hand, if $z\in Y\smallsetminus Y_0$, then by the construction of $Y$, there exists $j \leq k$ such that $z  \in \{x_{\gamma_j}, y _{\gamma_j}\}$.
        In particular, we have $ord(z) = \gamma_j$.
        Using $\gamma_j=ord(z) < \lambda_i$ together with the inequality (\ref{ineq for claim}), we conclude that $j \leq i$ and $\gamma_j \leq \gamma_i < \lambda_i$.
        It is now clear that $z \leq x_{\gamma_j} \leq x_{\gamma_i}$, as desired.

        Now, since $z \in Y$ by assumption and $x_{\gamma_i} \in Y$ by construction, $z \leq x_{\gamma_i}$ entails $z R x_{\gamma_i}$.
        Recall our definition of $\gamma_i$, in particular, the fact that $x_{\gamma_i} \in U_i \subseteq V(\varphi) \smallsetminus V(\psi)$.
        By our induction hypothesis, it follows $x_{\gamma_i} \in V{\restriction}(\phi)\smallsetminus V{\restriction}(\psi)$.
        Thus, $zR x_{\gamma_i}$ yields $\mathfrak{N}, z \not \models \phi \to \psi$, as desired.
    \ebullet

    Now, let us prove the claim for an arbitrary $\theta \in S(\varphi)$ such that $ipd(\theta)= i+1$.
    Note that, by the definitions of the formulas in the language of $\bipc$ and of $ipd$, such a formula $\theta$ must be built up, using only $\land$ and $\lor$, from formulas in $S_{\to}(\varphi) \cup S_{\gets}(\varphi)$ whose implicative degree is at most $i+1$.
    By what we proved above and by our induction hypothesis, the claim follows readily. \qed \vs

    Recall that we have $u \in X \smallsetminus V(\varphi)$ with $ord(u)= \alpha$, and that $\{x_\alpha, y_\alpha\} \subseteq Y_0 \subseteq Y$ by construction.
    Since $\M, u \not \models \varphi$, the above claim yields $\mathfrak{N}, u \not \models \varphi$. 
    Thus, $\Y$ refutes $\varphi$, as desired.
\end{proof}

\begin{Theorem} \label{thm fmp}
    The logic $LFC = \lc + \beta(\F_0) +\J(\F_1) + \J(\F_2) +\J(\F_3)$ has the FMP.
\end{Theorem}
\begin{proof}
    It suffices to show that if $\varphi \notin LFC$, then there exists $\A \in \V_{LFC}^{< \omega}$ such that $\A \not \models \varphi$.
    Accordingly, we suppose that $\varphi \notin LFC$, i.e., that $\V_{LFC} \not \models \varphi$.
    Since varieties are generated by their finitely generated members, there exists $\B \in \V_{LFC}$ such that $\B \not \models \varphi$ and $\B$ is $n$-generated, for some $n \in \omega$.
    If $\B$ is finite, we are done, so let us assume otherwise.
    By the Subdirect Representation Theorem (see, e.g., \cite[Thm.\ II.8.6]{Sanka2}), we can suppose without loss of generality that $\B$ is SI.
    By combining Theorem \ref{Thm poset structure} with Theorem \ref{thm finite submodel}, $\varphi$ is refuted by a finite comb or hcomb.
    Since the algebraic duals of the finite combs and hcombs are all contained in $\V_{LFC}$ by Proposition \ref{prop finite si are equal}, the result follows from the bi-Esakia duality.
\end{proof}

With the goal of this subsection now achieved, we can easily deduce that the logics $Log(FC)$ and $LFC$ coincide, and prove Theorem \ref{main thm}:\vs

\noindent \textit{Proof of Theorem} \ref{main thm}. By Proposition \ref{prop finite si are equal}, we know that $(\V_{Log(FC)})_{SI}^{<\omega}=(\V_{LFC})_{SI}^{<\omega}$.
Since both logics have the FMP ($Log(FC)$ by definition, and $LFC$ by Theorem \ref{thm fmp}), it follows that $\V_{Log(FC)}=\V_{LFC}$, i.e., that $Log(FC)=LFC$.
Thus, $Log(FC)$ is a finitely axiomatizable logic with the FMP, and is therefore decidable.
Consequently, the problem of determining if a finitely axiomatizable logic is contained in $Log(FC)$ is also decidable.
Since an extension $L$ of $\lc$ is locally tabular iff $L \nsubseteq Log(FC)$ (see Theorem \ref{Thm:locally-tabular-main}), we have established the decidability of the problem of determining if a finitely axiomatizable extension of $\lc$ is locally tabular.\qed

\vs

\noindent 
{\bf Acknowledgment}\ \ The first author was supported by the grant 2023.03419.BD from the Portuguese Foundation for Science and Technology (FCT).  
The second author was supported by the proyecto PID$2022$-$141529$NB-C$21$ de investigaci\'on financiado por MICIU/AEI/ 10.13039/501100 011033 y por FEDER, UE. He was also supported by the Research Group in Mathematical Logic, $2021$SGR$00348$ funded by the Agency for Management of University and Research Grants of the Government of Catalonia.

\newpage
\bibliographystyle{plain}

\end{document}